\documentclass[11pt,reqno]{amsart}
\usepackage{enumitem} 
\usepackage[margin=1.2in,includeheadfoot]{geometry}

\usepackage{amssymb}
\usepackage{amsaddr}
\usepackage{latexsym}
\usepackage{amsmath}
\usepackage{mathrsfs}

\usepackage{upgreek}
\usepackage{mathabx}
\usepackage{esint}
\usepackage{xcolor}
\usepackage{bbm}
\definecolor{citeblue}{RGB}{0,0,150}
\definecolor{citegreen}{RGB}{0,110,0}
\usepackage[colorlinks = true,
            linkcolor = citeblue,
            urlcolor  =  citeblue,
            citecolor =  citegreen,
        ]{hyperref}   
\usepackage{empheq}
\usepackage{cleveref}
\usepackage{mathtools}
\usepackage{amsthm}
\usepackage[title]{appendix}
\usepackage{amsfonts}
\usepackage{relsize}
\usepackage{esvect}

\makeatletter
\def\l@subsection{\@tocline{2}{0pt}{2pc}{0pc}{}}
\makeatletter
\usepackage{cases}
\makeatletter

\newcommand{\Rmn}[1]{\expandafter\@slowromancap\romannumeral #1@}
\makeatother

\usepackage{graphicx}
\usepackage[left,modulo]{lineno}
\newtheorem{thm}{Theorem}[section]

 \newtheorem{Th}[thm]{Theorem}

 \newtheorem{Prop}[thm]{Proposition}
 \newtheorem{Co}[thm]{Corollary}
\newtheorem{Lm}[thm]{Lemma}

\newtheorem{Dfi}[thm]{Definition}
\newtheorem{Rm}[thm]{Remark}

\numberwithin{equation}{section}
\newcommand{\be}{\begin{equation}}
\newcommand{\ee}{\end{equation}}
\newcommand{\bg}{\begin{gather}}
\newcommand{\eg}{\end{gather}}
\newcommand{\ba}{\begin{align}}
\newcommand{\ea}{\end{align}}
\newcommand{\bad}{\begin{aligned}}
\newcommand{\ead}{\end{aligned}}
\newcommand{\R}{\mathbb{R}}
\newcommand{\N}{\mathbb{N}}

\newcommand{\mca}[1]{\mathcal{#1}}

\newcommand{\nf}{\infty}

\def\avint{\mathop{\mathchoice{\,\rlap{-}\!\!\int}
                              {\rlap{\raise.15em{\scriptstyle -}}\kern-.2em\int}
                              {\rlap{\raise.09em{\scriptscriptstyle -}}\!\int}
                              {\rlap{-}\!\int}}\nolimits}

\def\avint{\mathop{\,\rlap{-}\!\!\int}\nolimits}

\def\smi{\setminus}

\def\XXint#1#2#3{{\setbox0=\hbox{$#1{#2#3}{\int}$ }
\vcenter{\hbox{$#2#3$ }}\kern-.6\wd0}}

\def\wih{\widehat}
\usepackage{dsfont}
\def\mkt{\mkern2mu}
\def\mko{\mkern1mu}

\def\ga{\gamma}
\def\La{\Lambda}

\def\lf{\left}

\def\de{\delta}
\def\rg{\right}

\def\al{\alpha}

\def\la{\lambda}

\def\wti{\widetilde}

\def\vae{\varepsilon}

\def\vp{\varphi}

\def\Om{\Omega}

\def\p{\partial}

\def\vt{\vartheta}
\def\bg{\vec{g}}

\def\Ga{\Gamma}

\def\vp{\varphi}

\def\si{\sigma}

\newcommand{\metricsub}[3]{_{
    \mkern-#3mu
    \smash{\lower #1\hbox{$\scriptscriptstyle #2$}}
}}

\newcommand{\ov}[1]{\overline{#1}}

\DeclareMathOperator{\diam}{diam}
\DeclareMathOperator{\dist}{dist}

\newcommand{\lap}{\Delta}

\begin{document}
\title[Uniform estimates with critical Sobolev coefficients]{Uniform estimates for elliptic systems with critical Sobolev coefficients}
\author{\small Tian Lan}
\address{\smaller Department of Mathematics, ETH Zürich, 101 Rämistrasse, 8092 Zürich, Switzerland}
\date{\today}
\begin{abstract}
The elliptic regularity theory with $W^{1,n}$ coefficients dates back to Miranda, who established $W^{2,2}$ estimates for elliptic equations in divergence and non-divergence forms. 
For any such fixed operator, the $L^p$ estimates for all $1<p<\infty$ follow from the theory of VMO coefficients developed since the 1990s.
However, the VMO theory does not yield constants in the estimates uniformly controlled by the $W^{1,n}$ norm of the coefficients, since the total $W^{1,n}$ norm does not control the VMO modulus. 
In this work, we establish uniform estimates for divergence-form elliptic systems with coefficients in the critical Sobolev--Lorentz spaces $W^{1,(n,s)}$ with $1\le s<\infty$, allowing an unbounded spatially skew-symmetric part.
The constants in the estimates depend quantitatively on the critical Sobolev norms, and we show that they are sharp for scalar equations when $p>n$ and for systems when $p\neq 2$. 
\end{abstract} 
\maketitle
 \vskip0.5cm
\tableofcontents

\section{Introduction}\label{sec:intro}
For simplicity, we first consider scalar equations in divergence form. 
Let $\Om\subset \R^n$ be a bounded domain with sufficiently smooth boundary, and let $A\coloneq\{A^{\al\beta}\}$ be a bounded measurable $n\times n$ matrix of coefficients satisfying a uniform ellipticity condition. 
We consider the following Dirichlet problem for the operator $L_Au=\p_\al(A^{\al\beta}\p_\beta u)$:
    \begin{equation}\label{eqdiv}
    \begin{dcases}
    \begin{aligned}
       L_Au&\,=\, \p_\al f_\al&& \text{ in }\;\Om,\\
        u&\,=\,0&& \text{ on }\,\p \Om.
        \end{aligned}
    \end{dcases}
\end{equation}
For $1<p<\nf$ and $\{f_\al\}_{1\le\al\le n}\subset L^p(\Om)$, we are interested in existence and uniqueness of $u\in W^{1,p}_0(\Om)$ solving~\eqref{eqdiv} and validity of the following estimate:
\begin{equation}\label{scalestW1p}
    \|Du\|_{L^p(\Om)}\le C \sum_{\al=1}^n\|f_\al\|_{L^p(\Om)}.
\end{equation}
This is equivalent to 
\begin{equation}\label{ine:invLAbd}
   \|u\|_{W^{1,p}(\Om)}\le C\|L_A u\|_{W^{-1,p}(\Om)}, \qquad \forall\, u\in W^{1,p}_0(\Om),
\end{equation}
where $C>0$ is a constant depending only on $\Om$, $p$, and $A$.
When $A$ is uniformly continuous, the estimate~\eqref{ine:invLAbd} holds for every $p\in (1,\nf)$; see~\cite{CalZyg52,Morrey08,Simader72,Campanato80}. 
For merely bounded measurable coefficients, Meyers~\cite{Meyers63} proved that there exists $\vae=\vae(\Om,A)>0$ such that, for every $p\in (2-\vae,2+\vae)$, the problem~\eqref{eqdiv} admits a unique solution in $W^{1,p}_0(\Om)$ satisfying~\eqref{scalestW1p}.

A natural class of discontinuous coefficients is given by the critical Sobolev space $W^{1,n}$. 
Elliptic regularity for such coefficients goes back to Miranda~\cite{Miranda63}, who established $W^{2,2}$ estimates for elliptic equations in divergence and non-divergence forms with $W^{1,n}$ coefficients in dimensions $n\ge3$.
Alvino and Trombetti~\cite{AlviTromb84} subsequently obtained analogous estimates under a smallness assumption on the $L^{n,\nf}$ norms (see Definition~\ref{def-Lor}) of the coefficient derivatives. 
Along the same line, Chiarenza and Franciosi~\cite{ChiaFran92} generalized Miranda's theorem under a vanishing Morrey condition on the coefficient derivatives.
Based on Meyers' work~\cite{Meyers63}, $L^p$ estimates for divergence and non-divergence form equations with $W^{1,n}$ coefficients were also obtained in~\cite{MandPorr90} for $p$ sufficiently close to $2$.

A broader framework is provided by coefficients of vanishing mean oscillation (VMO), for which elliptic estimates are available throughout the range $1<p<\nf$.
Chiarenza, Frasca, and Longo developed the $W^{2,p}$ theory for equations in non-divergence form with VMO coefficients, while Di Fazio~\cite{DiFazio96}, Auscher and Qafsaoui~\cite{AusQaf02} established the corresponding $W^{1,p}$ estimates in divergence form.

In 1998, Caffarelli and Peral~\cite{CaffaPeral98} introduced an approximation method using a real-variable argument from the Calderón--Zygmund decomposition to obtain $W^{1,p}$ estimates for divergence-form equations. 
Building on this approach, Byun~\cite{Byun05} and Byun--Wang~\cite{ByunWang04} obtained $W^{1,p}$ estimates in Lipschitz domains and, more generally, Reifenberg domains under smallness assumptions on the boundary flatness and the BMO seminorms of the coefficients.
Analogous estimates for elliptic systems satisfying the Legendre condition were obtained in~\cite{ByunChenKimWang07,ByunWang08}. 
For higher-order systems in divergence and non-divergence forms under the Legendre--Hadamard condition, an $L^p$ solvability theory was developed by Dong and Kim~\cite{DongKim11A}.
Related results are also available for partially VMO or partially BMO coefficients; see~\cite{Krylov09,DongKim10,ByunRyuWang10,DongKim11B}.

For any ball $B\subset \R^n$ with $n\ge 2$ and $f\in W^{1,(n,\nf)}(B)$, the Poincar\'e inequality gives
\begin{equation}\label{ine:Poincare}
    \fint_B \big|f-(f)_B\big|^2\le C\|Df\|_{L^{n,\nf}(B)}^2,
\end{equation}
where $C=C(n)$. Thus, on a bounded domain $\Om$ with sufficiently smooth boundary, the mean-oscillation hypotheses in the aforementioned elliptic estimates are satisfied for a coefficient matrix $A=\{A^{\al\beta}\}$ whenever $\|DA\|_{L^{n,\nf}(\Om)}$ is sufficiently small or $\|DA\|_{L^{n,\nf}(\Om\cap B_r(x))}$ tends to zero uniformly in $x\in\bar\Om$ as $r\to 0$.
The latter condition holds for any fixed coefficient matrix $A\in W^{1,n}(\Om)$ and more generally for $A\in W^{1,(n,s)}(\Om)$ with $1\le s<\nf$.
For recent applications of elliptic estimates involving critical Sobolev--Lorentz coefficients, see for instance~\cite{LaManLeoSch20,MoscaPasca24,MarRiv25,Pascale26,BerLanMarRiv26,Lan26}.

We now consider a sequence of operators $\{L_{A_k}\}_{k\in\N}$ given by $L_{A_k}u=\p_\al(A_k^{\al\beta}\p_\beta u)$ on a fixed bounded domain $\Om$ with sufficiently smooth boundary.
We assume that all $A_k=\{A_k^{\al\beta}\}$ satisfy a uniform ellipticity condition and that $\|DA_k\|_{L^n(\Om)}\le\eta$ for some constant $\eta\ge 0$.
A natural question is whether the following estimate holds:
\begin{equation}\label{ine:toyuniest}
   \|u\|_{W^{1,p}(\Om)}\le C\|L_{A_k} u\|_{W^{-1,p}(\Om)}, \qquad \forall\, u\in W^{1,p}_0(\Om),
\end{equation}
with $C$ independent of $k$ and $1<p<\nf$.
We note that the elliptic theory for BMO or VMO coefficients discussed above does not directly imply~\eqref{ine:toyuniest} for an arbitrary $\eta>0$.
The constants in the VMO theory depend on the rate of decay of the mean oscillation of $A_k$, whereas the bound on $\|DA_k\|_{L^n(\Om)}$ alone does not control this rate uniformly in $k$.
A typical example is obtained as follows: for $k\ge 1$, set $\Om=B_1$ and $A_k(x)=A_0(kx)$, where $A_0\in L^\nf(\R^n,\R^{n\times n})$ is a symmetric uniformly elliptic matrix of coefficients with $\|DA_0\|_{L^n(\R^n)}\le\eta$.
A change of variables gives $\|DA_k\|_{L^n(B_1)}=\|DA_0\|_{L^n(B_{k})}\le\eta$. 
However, we have
\[
    \fint_{B_{\frac 1k}}\Big|A_k-(A_k)_{B_{\frac 1k}}\Big|^2
    =
    \fint_{B_1}\Big|A_0-(A_0)_{B_1}\Big|^2>0.
\]
Hence the sequence $\{A_k\}$ does not admit a uniform VMO modulus, and we cannot obtain~\eqref{ine:toyuniest} by directly applying the VMO theory. Nevertheless, for this toy model, we have
\[
    \int_{B_1\smi B_{1/4}} |DA_k|^n\to 0\quad \text{as }\;k\to\nf.
\]
Thus, the small BMO theory yields uniform $L^p$ gradient estimates on $B_1\smi B_{1/2}$. The corresponding bounds on $B_{1/2}$ then follow from the whole-space estimate in~\cite{IwaSbo98} by scaling and a cut-off argument. This proves~\eqref{ine:toyuniest} for this family. 
However, the general case allows concentration at multiple points and at nested scales, with microscopic structures that may vary along the sequence. 

We would like to compare the problem~\eqref{ine:toyuniest} with the elliptic homogenization problems, where one seeks estimates independent of $\vae>0$ for coefficients of the form $A_\vae(x)=A(x/\vae)$, with $A$ fixed periodic or almost-periodic; see~\cite{AvellaLin87,AvellaLin91,CaffaPeral98,ArmShen16,Shen17,Shen18}.
Unless $A$ is constant, the family $\{A_\vae\}$ has no common VMO modulus, and $\|DA_\vae\|_{L^n(B_1)}\to\nf$ as $\vae\to0$ even if $A$ is smooth.  
Nevertheless, the common microscopic structure provides uniform estimates at microscopic scales, which serve as fundamental ingredients for establishing the global estimates. 
In our problem, the critical Sobolev seminorm of coefficients is uniformly bounded, but no common microscopic structure is prescribed. We establish an analogue of~\eqref{ine:toyuniest} for elliptic systems whose coefficients have a bounded spatially symmetric part and whose derivatives have $L^{n,s}$ norms bounded by $\eta$, with constants depending quantitatively on $\eta$.

As discussed above, to prove such a uniform estimate we need to exploit properties of the $W^{1,n}$ space beyond its embedding into VMO.
This phenomenon also appears in the application of the John--Nirenberg inequality~\cite{JohnNiren61}.
For $f\in L^1(B_1) \cap \mathrm{BMO}(B_1)$, the John--Nirenberg inequality gives constants $\ga=\ga(n)>0$ and $C=C(n)>0$ such that for any ball $B\subset B_1$,
\begin{equation}\label{JohnNirine}
    \fint_B \exp\lf(\ga \frac{|f- (f)_B|}{[f]_{\mathrm{BMO}(B)}}\rg)\le C.
\end{equation}
In particular, for $f\in W^{1,n}(B_1)\subset\mathrm{VMO}(B_1)$ this implies $e^f\in L^1(B_1)$. 
However, the bound obtained in this way depends on the VMO modulus of $f$, which is not fully controlled by $\|f\|_{W^{1,n}(B_1)}$.
In contrast, the Moser--Trudinger inequality~\cite{Trudinger67,Moser71}, combined with Young's inequality, gives a uniform estimate of $\|e^f\|_{L^1(B_1)}$ in terms of $\|f\|_{W^{1,n}(B_1)}$.
An analogous exponential integrability estimate also appears in our proof; see Lemma~\ref{lm:Wolff-exp}.

Our proof follows the real-variable approach developed by Shen in~\cite{Shen05,Shen07}, which originated in the work of Caffarelli--Peral~\cite{CaffaPeral98}. This reduces the $L^p$ gradient estimates to local reverse H\"older inequalities for homogeneous solutions.
To obtain these inequalities with uniform constants, we adapt the gradient estimate method developed by Luc Nguyen~\cite{Nguyen23} and inspired by Li~\cite{Li17},
and control the coefficient term pointwise by a Wolff potential without imposing a Dini mean oscillation condition. 

\subsubsection*{Acknowledgements.}
The author thanks Prof. Tristan Rivi\`ere for his support. 
The author also thanks Prof. Luc Nguyen for helpful correspondence.

\subsubsection*{AI statement.}
DeepSeek V4 Pro was used for preliminary literature searches on gradient estimates for homogeneous solutions, Hajłasz spaces, and pathological solutions, including~\cite{Li17,Nguyen23,KoskYangZhou11,Zhou11,JinMazyaVanSch09}.
DeepSeek V4 Pro, ChatGPT-5.6 Sol, and ChatGPT-6.0 Sol were used as auxiliary editorial tools for refining the exposition and correcting typos in Sections~\ref{sec:intro} and~\ref{sec:mainres}. 
After the proofs of Theorems~\ref{th:divest},~\ref{Th:divHajespa}, and Proposition~\ref{prop:shasca} were completed, ChatGPT-6.0 Sol was used to assist in the algebraic computation of the coefficient tensor $A^{\al\beta}_{ij}$ in~\eqref{exsysA}; the validity of computation was verified by the author.
Apart from these uses, all the mathematical ideas, arguments, and proofs are developed entirely by the author.
\section{Main results}\label{sec:mainres}
We first introduce some notation and definitions used in the paper. Let $n\ge 2$. We write $B_r=B_r(0)$. For $x\in \R^n$, we denote $x$ by $(x',x_n)$ where $x'\in \R^{n-1}$. For a measurable set $U\subset \R^n$ with $0<|U|<\nf$ and $f\in L^1(U)$, we write
\[
    (f)_U\coloneq \fint_U f= \frac{1}{|U|}\int_Uf.
\]
We recall the definitions of the function and distribution spaces used throughout this work.
\begin{Dfi}\label{def-Lor}
Let $\Om\subset\R^{n}$ be a bounded open set. Given a measurable function $f\colon \Om\to\R$, we define the distribution function and the decreasing rearrangement of $f$ as
\begin{equation}\label{disdecreadef}
    d_{f}(\la)\coloneq\mca{L}^{n}\bigl\{x\in \Om:|f(x)|>\la\bigr\},
    \qquad
    f^{*}(t)\coloneq\inf\bigl\{\la\ge 0 : d_{f}(\la)\le t\bigr\}.
\end{equation}
For $0< p<\nf$ and $0< q\le \infty$, we define the Lorentz quasinorm
\begin{equation}\label{eq:defLor}
\|f\|_{L^{p,q}(\Om)}
\coloneq\bigl\|t^{\frac1p}f^{*}(t)\bigr\|_{L^{q}(\R_{+},\,dt/t)}
    =p^{\frac1q}\bigl\|\la\,d_{f}(\la)^{\frac1p}\bigr\|_{L^{q}(\R_{+},\,d\la/\la)}.
\end{equation}
The Lorentz space $L^{p,q}(\Om)$ consists of all measurable $f$ with
$\|f\|_{L^{p,q}(\Om)}<\infty$. In particular, we have $L^{p,p}(\Om)=L^p(\Om)$.
When $p>1$ and $1\le q\le \nf$, the Lorentz quasinorm is equivalent to a norm; see for instance \cite[Ch.~4]{BennShar88}. For $1\le s<\nf$, we define the Sobolev--Lorentz space
\begin{equation}\label{eq:defSoblor}
    W^{1,(n,s)}(\Om)\coloneq \{f\in L^{n,s}(\Om): Df\in L^{n,s}(\Om)\}.
\end{equation}
For $1<p<\nf$, we set $p'=p/(p-1)$, and define the negative-order Sobolev space $W^{-1,p}(\Om)$ as the dual space of $W^{1,p'}_0(\Om)$ with $\|f\|_{W^{1,p'}_0(\Om)}\coloneq \|Df\|_{L^{p'}(\Om)}$.
\end{Dfi}
Our main assumption on the domain is that it is locally a graph of a Lipschitz continuous function with a small Lipschitz constant.
\begin{Dfi}\label{def:smaLip}
    Let $\Om\subset \R^n$ be a bounded open set and let $R_0>0$. We say $\Om$ is $(\de, R_0)$-Lipschitz if for any $x_0\in \p\Om$, there exists a Lipschitz continuous function $\ga:\R^{n-1}\to \R$ such that, upon rotating and relabeling the coordinate axes if necessary,
    \[
        \Om\cap B_{R_0}(x_0)=\{x\in B_{R_0}(x_0):x_n>\ga(x') \},
    \]
    and $|\ga(x')-\ga(y')|\le \de |x'-y'|$ for any $x',y'\in \R^{n-1}$ with $\max\big\{|x'-x_0'|,|y'-x_0'|\big\}\le R_0$.
\end{Dfi}

Let $n\ge 2$, $m\ge 1$, $1\le s<\nf$, and let $\Om\subset \R^n$ be a bounded Lipschitz domain. Suppose $A=\{A^{\al\beta}_{ij}\}$ with $A^{\al\beta}_{ij}\in W^{1,(n,s)}(\Om,\R)$ for all $1\le i,j\le m$ and $1\le \al,\beta\le n$. We write
\begin{equation}\label{def:spatparts}
 S^{\al\beta}_{ij}\coloneq\frac{A^{\al\beta}_{ij}+A^{\beta\al}_{ij}}2,
 \qquad K^{\al\beta}_{ij}\coloneq\frac{A^{\al\beta}_{ij}-A^{\beta\al}_{ij}}2.
\end{equation} 
Let $\La\ge 1$ and $\eta\ge 0$ be constants. In this paper, we assume 
\begin{equation}\label{DAbdLns}
   \big\|DA^{\al\beta}_{ij}\big\|_{L^{n,s}(\Om)}\le\eta,\qquad  \forall\,1\le i,j\le m,\;\; \forall \,1\le\al,\beta\le n.
\end{equation}
We also assume the spatially symmetric part of $A$ is bounded:
\begin{equation}\label{bdsympart}
    \big\|S^{\al\beta}_{ij}\big\|_{L^\nf(\Om)}\le\La,\qquad  \forall\,1\le i,j\le m,\;\; \forall \,1\le\al,\beta\le n.
\end{equation}
For $\vp_1\in C^\nf (\ov\Om,\R^m)$ and $\vp_2\in C_c^\nf(\Om, \R^m)$, we define
\begin{equation}\label{eq:biliform}
 \mca B_A(\vp_1,\vp_2)\coloneq\int_\Om A^{\al\beta}_{ij}\p_\beta \vp_1^j\mko\p_\al\vp_2^i,
\end{equation}
where $\vp_k=(\vp_k^i)_{1\le i\le m}$ for $k=1,2$. By~\eqref{bdsympart}, for any $1<p<\nf$, we have
\[
   \bigg| \int_\Om S^{\al\beta}_{ij}\p_\beta \vp_1^j\mko\p_\al\vp_2^i\bigg|\le C\mko\La \mko\|\vp_1\|_{W^{1,p}(\Om)} \|\vp_2\|_{W^{1,p'}_0(\Om)},
\]
where $C=C(n,m)$. 
In addition, since $W^{1,(n,s)}(\R^n)\subset \mathrm{BMO}(\R^n)$, the div--curl lemma~\cite{CoiLioMeySem93}, combined with the Sobolev extensions, implies that 
\begin{align*}
     \bigg| \int_\Om K^{\al\beta}_{ij}\p_\beta \vp_1^j\mko\p_\al\vp_2^i\bigg|&=\frac 12\bigg| \int_\Om K^{\al\beta}_{ij}\big(\p_\beta \vp_1^j\mko\p_\al\vp_2^i-\p_\al \vp_1^j\mko\p_\beta\vp_2^i\big)\bigg|\\
     &=\frac 12\bigg| \int_\Om \Big(K^{\al\beta}_{ij}-\fint_\Om K^{\al\beta}_{ij}\Big)\big(\p_\beta \vp_1^j\mko\p_\al\vp_2^i-\p_\al \vp_1^j\mko\p_\beta\vp_2^i\big)\bigg|\\
     &\le C \mko\eta \mko\|\vp_1\|_{W^{1,p}(\Om)} \|\vp_2\|_{W^{1,p'}_0(\Om)},
\end{align*}
where $C=C(\Om,m,p)$ and the second equality follows from integration by parts. 
It follows that 
\begin{equation}\label{bdBALp}
    \big|  \mca B_A(\vp_1,\vp_2)\big|\le C\mko (\La+\eta)\|D\vp_1\|_{L^{p}(\Om)} \|\vp_2\|_{W^{1,p'}_0(\Om)}.
\end{equation}
Moreover, the condition $S^{\al\beta}_{ij}\in L^\nf(\Om)$ is in fact necessary for the estimate 
\[
    \big|  \mca B_A(\vp_1,\vp_2)\big|\le C\|\vp_1\|_{W^{1,2}(\Om)} \|\vp_2\|_{W^{1,2}_0(\Om)}
\]
to hold. Indeed, fixing $1\le i,j\le m$ and setting $\vp_1=u\mko e_i$, $\vp_2=v\mko  e_j$, where $(e_i)_{i=1}^m$ is the standard basis of $\R^m$, reduces the inequality to the case of scalar functions. 
The necessity of boundedness of the symmetric part then follows from the proof of~\cite[Prop.~3.1]{MazyaVerb06}. 
See also~\cite{LiPipher19,HofLiMayPi22} for related settings of the coefficients for scalar elliptic equations.

Now for $u=(u^j)_{1\le j\le m}$, we define $(L_A u)^i=\p_{\al}\big( A^{\al\beta}_{ij}\p_\beta u^j\big)$.
The inequality~\eqref{bdBALp} is then equivalent to
\begin{equation}\label{eq:W-1pestLA}
    \|L_A u\|_{W^{-1,p}(\Om)}\le C(\La+\eta) \|Du\|_{L^{p}(\Om)}.
\end{equation}
We assume an integral ellipticity condition on $A$, which can be written as
\begin{equation}\label{eq:intcond}
 \mca B_A(\vp,\vp)\ge\La^{-1}\|D\vp\|_{L^2(\Om)}^2,
 \qquad \forall \,\vp\in W^{1,2}_0(\Om,\R^m).
\end{equation}
The well-definedness of~\eqref{eq:intcond} follows from~\eqref{bdBALp}. This integral ellipticity condition is weaker than the Legendre condition, but is stronger than the Legendre--Hadamard condition: for almost every $x\in \Om$, there holds
\begin{equation}\label{eq:LegHadcon}
    A^{\al\beta}_{ij}(x)\eta^i\eta^j \xi_\al \xi_\beta\ge \La^{-1}|\eta|^2|\xi|^2,\qquad \forall\,\eta=(\eta^i)\in \R^m,\;\;\forall\,\xi=(\xi_\al)\in \R^n.
\end{equation}
See for instance the proof of~\cite[Thm.~1.1]{Dindos21}. We now state our main result.

\begin{Th}\label{th:divest}
Let $1<p<\nf$. There exists $\de=\delta(n,m,\La,p)>0$ such that the following holds. Suppose $\Om\subset\R^n$ is bounded and $(\de,R_0)$-Lipschitz for some $R_0>0$, and suppose $A=\{A^{\al\beta}_{ij}\}$ satisfies~\eqref{DAbdLns},~\eqref{bdsympart}, and~\eqref{eq:intcond}.
Then for any $f=\{f_i^\al\}\in L^{p}(\Om,\R^{m\times n})$, there exists a unique $u=(u^j)\in W^{1,p}_0(\Om,\R^m)$ solving
\begin{equation}\label{sysinth}
    \p_{\al}\big( A^{\al\beta}_{ij}\p_\beta u^j\big)\,=\,\p_\al f^\al_i\qquad\text{in }\Om.
\end{equation}
Moreover, there exist positive constants $C=C(\Om)$, $C_0=C_0(\La,n,m)$, and $C_1=C_1(\La,n,m,s)$ such that
\begin{equation}\label{W1pestW1n}
 \|Du\|_{L^p(\Om)}\le C\exp\Big(C_0 \max(p,p') +C_1\max (p,p')^{s}\eta^s\Big)\mko\|f\|_{L^{p}(\Om)}.
\end{equation}
\end{Th}
\begin{Rm}
    For $s>1$, the weak type $(1,1)$ estimate fails in general for uniformly elliptic divergence-form operators with coefficients in $W^{1,(n,s)}$, in contrast to the standard Laplacian and, more generally, to operators whose coefficients satisfy
    \[
       \fint_{\Om\cap B_r(x)} \big|A-(A)_{\Om\cap B_r(x)}\big|\le C\log^{-2} (r),\qquad \forall\, 0<r<\frac12,\;\forall\, x\in \ov\Om,
        \]
 see~\cite{DongKim17,DongEscaKim18}. Indeed, the symmetric coefficient matrix $A=\{A^{\al\beta}\}$ constructed in~\cite{JinMazyaVanSch09} belongs to $W^{1,(n,s)}(B_1)$ for any $s>1$. 
 Suppose that for any $f=\{f_\al\}\in L^2(B_1,\R^n)$ and $u\in W^{1,2}_0(B_1)$ solving
    \begin{equation}\label{examuf}
      \p_\al(A^{\al\beta}\p_\beta u)=\p_\al f_\al \quad \text{in }B_1,
\end{equation}
the following estimate holds:
\[
    \|Du\|_{L^{1,\nf}(B_1)}\le C \|f\|_{L^1(B_1)}.
\]
Then the Marcinkiewicz interpolation and duality imply that for any $p>2$ and $f\in L^p(B_1,\R^n)$, the solution $u\in W^{1,2}_0(B_1)$ of~\eqref{examuf} satisfies
\begin{equation}\label{Lpmarcinassu}
    \|Du\|_{L^p(B_1)}\le C p\|f\|_{L^p(B_1)},
\end{equation}
where $C$ is independent of $p$ and $f$.
However, multiplying the solution constructed in~\cite[Prop.~1.6]{JinMazyaVanSch09} by a cut-off function gives a function $u_0\in \bigcap_{1<p<\nf}W^{1,p}_0(B_1)$ satisfying 
\begin{equation*}
\begin{dcases}
    \p_\al(A^{\al\beta}\p_\beta u_0)\in C_c^\nf (B_1),\\
    \|Du_0\|_{L^p(B_1)}\ge c\mko p^2,\quad \forall\, 2<p<\nf,
    \end{dcases}
\end{equation*}
where $c>0$ is independent of $p$. This contradicts the estimate~\eqref{Lpmarcinassu} and disproves the weak type $(1,1)$ estimate for $s>1$. 

To the author's knowledge, it remains open whether the weak type $(1,1)$ estimate holds when $s=1$ and $\Om$ is smooth.
However, under the assumptions of Theorem~\ref{th:divest} with $s=1$ and $\Om$ having a $C^{1,\textup{Dini}}$ boundary (see e.g.~\cite{DongEscaKim18}), adapting the proofs of Lemmas~\ref{lm:Nguint} and~\ref{lm:Ngubdy} to such domains, together with~\eqref{eq:Wolend} and duality, yields
\begin{equation}\label{Lnn-1nfbd}
    \|u\|_{L^{n/(n-1),\nf}(\Om)}\le C e^{ C\eta} \|f\|_{L^1(\Om)},
\end{equation}
where $C=C(n,m,\La,\Om)>0$.
Motivated by this, we conjecture that the weak type $(1,1)$ estimate holds when $\Om$ has a $C^{1,\textup{Dini}}$ boundary and $s=1$, with a quantitative bound similar to~\eqref{W1pestW1n} and~\eqref{Lnn-1nfbd}.
\end{Rm}
A natural corollary of Theorem~\ref{th:divest} is the following uniform estimate with coefficients in higher-order critical Sobolev spaces, which follows from the same argument as in~\cite[Lem.~III.1]{Lan26} by using the Sobolev--Lorentz embeddings and product estimates (see e.g.~\cite[Lem.~II.6]{Lan26}).
\begin{Co}
    Let $k\in \N^+$ with $k\le n-1$, $1< p<\frac nk$, $1\le s<\nf$, and let $\Om\subset\R^n$ be a bounded $C^{k,1}$ domain. 
    Suppose $A=\{A^{\al\beta}_{ij}\}$ satisfies~\eqref{bdsympart},~\eqref{eq:intcond}, and 
    \[
       \sum_{\ell=1}^k  \|D^{\ell} A\|_{L^{ n/\ell,s}(\Om)}\le \eta.
    \]
    Then for any $F=(F^i)\in W^{k-1,p}(\Om,\R^m)$, there exists a unique $u=(u^j)\in W^{k+1,p}(\Om,\R^m)$ solving
\[
    \begin{dcases}
    \begin{aligned}
      \p_{\al}\big( A^{\al\beta}_{ij}\p_\beta u^j\big)&\,=\,  F^i&& \text{ in }\;\Om,\\
        u^j&\,=\,0&& \text{ on }\,\p \Om.
        \end{aligned}
    \end{dcases}
\]
Moreover, there exists a constant $C=C(\La,\Om,m,p,s)>0$ such that
\begin{equation*}
 \|u\|_{W^{k+1,p}(\Om)}\le C\exp(C\eta^s)\mko\|F\|_{W^{k-1,p}(\Om)}.
\end{equation*}
\end{Co}
Inspired by discussions of non-divergence form elliptic equations with coefficients in the critical \textit{Aronszajn--Gagliardo--Slobodeckij space} $W^{\theta,n/\theta}$~\cite{ChiFraLon91,ChiFraLon93}, we extend our uniform estimates for divergence-form operators to coefficients in these fractional Sobolev spaces, and more generally, in critical \textit{fractional Hajłasz--Lorentz spaces}. 

Let $0<\theta<1$, $\Om$ be a Lipschitz domain, and we set $X=L^{n/\theta,s}(\Om)$. 
Following the notation of~\cite[Sec.~7.1]{Heikkinen18}, the Hajłasz--Lorentz space $M^{\theta,X}(\Om)$ consists of all $A\in X$ such that there exist a nonnegative function $\wti g\in X$ and a set $E\subset \Om$ with $|E|=0$ such that
\begin{equation}\label{conwtig}
    |A(x)-A(y)|\le |x-y|^\theta \big(\wti g(x)+\wti g(y)\mko \big),\qquad \forall\, x,y\in \Om\smi E.
\end{equation} 
In particular, when $s=n/\theta$, this coincides with the Triebel--Lizorkin space $F^{\theta}_{n/\theta, \nf}(\Om)$ and contains $W^{\theta,n/\theta}(\Om)$; see~\cite{KoskYangZhou11,Rychkov99,Triebel10}.
We now state our uniform estimate with coefficients in $M^{\theta,X}(\Om)$.
\begin{Th}\label{Th:divHajespa}
    Let $1<p<\nf$, and let $\de$, $\Om$ be as in Theorem~\ref{th:divest}. 
   Suppose $A=\{A^{\al\beta}_{ij}\} $ satisfies~\eqref{bdsympart} and~\eqref{eq:intcond}. 
   Assume further that there exist a nonnegative function $\wti g\in L^{n/\theta,s}(\Om)$ and a set $E\subset \Om$ with $|E|=0$ such that~\eqref{conwtig} holds and
   \begin{equation}\label{triziass}
       \|\wti g\|_{L^{n/\theta,s}(\Om)}\le \eta.
   \end{equation}
Then for any $f=\{f_i^\al\}\in L^{p}(\Om,\R^{m\times n})$, there exists a unique $u=(u^j)\in W^{1,p}_0(\Om,\R^m)$ solving~\eqref{sysinth}, and there exist positive constants $C=C(\Om)$, $C_2=C_2(\La,n,m)$, and $C_3=C_3(\La,n,m,\theta,s)$ such that
\begin{equation}\label{W1pestHaj}
 \|Du\|_{L^p(\Om)}\le C\exp\Big(C_2 \max(p,p') +C_3\max (p,p')^{s}\eta^{s}\Big)\mko\|f\|_{L^{p}(\Om)}.
\end{equation}
\end{Th}
In Sections~\ref{sec:intest} and~\ref{sec:global}, we focus on the proof of Theorem~\ref{th:divest}. 
Using a very similar argument, we provide a brief proof of Theorem~\ref{Th:divHajespa} at the end of Section~\ref{sec:global}. In the last section, we prove that the quantitative dependence on $\eta$ in~\eqref{W1pestW1n} is sharp.
\section{Interior estimates}\label{sec:intest}
In this section, we prove the reverse H\"older inequality on interior balls required for proving~\eqref{W1pestW1n}. Before establishing the elliptic estimates, we first record an estimate for the Wolff potentials.
\subsection{On Wolff potentials}
\
\vskip5pt
Let $n\ge 2$, $p>1$, and $0<\al<n/p$. For a nonnegative function $G\in L^1(\R^n)$ with compact support, we define the \textit{Wolff potential} by
\begin{equation}\label{def:Wolffpot}
    \mathbf W_{\al,p}G(x)\coloneq\int_0^\nf\bigg(t^{\al p}\fint_{B_t(x)}G\bigg)^{\frac1{p-1}}\,\frac{dt}{t}.
\end{equation}
The \textit{Riesz potential} is given by
\[
 I_\al G(x)\coloneq \frac{1}{c_{n,\al}}\int_{\R^n}G(y)\mko |x-y|^{\al-n}\,dy,
\]
where $c_{n,\al}$ is a positive constant. Since $\al p<n$, the classical pointwise comparison between the Wolff and Havin--Maz'ya potentials gives
\begin{equation}\label{comWolHaMaz}
    \mathbf W_{\al,p}G
    \le C I_\al\big(\mko(I_\al G)^{\frac1{p-1}}\big),
\end{equation}
where $C=C(n,\al,p)$; see~\cite[Part II, Lem.~7.1]{Mazya72} or~\cite[Thm.~3.1]{AdaMey72}.
We record the following exponential integrability estimate for Wolff potentials in the critical case; see also~\cite[Rem.~3.2 and Thm.~A]{CianSchw18}.
\begin{Lm}\label{lm:Wolff-exp}
Let $p>1$, $0<\al<n/p$, $1\le s<\nf$, and set
\[
    q\coloneq\frac{n(p-1)}{\al p}.
\]
Let $\eta\ge0$ and let $B\subset\R^n$ be a ball. Let $g\ge0$ be supported in $B$ with $\|g\|_{L^{q,s}(B)}\le\eta$.
Then there exists a constant $C_4=C_4(n,\al,p,s)>0$ such that for any $\ga>0$,
\begin{equation}\label{eq:WolLor}
    \fint_B\exp\big(\ga\mathbf W_{\al,p}(g^{p-1})\big)
    \le 4\exp\big(C_4(\ga\eta)^s\big).
\end{equation}
In addition, when $s=1$, we have $\mathbf W_{\al,p}(g^{p-1})\in C^0(\ov B)$ and
\begin{equation}\label{eq:Wolend}
    \|\mathbf W_{\al,p}(g^{p-1})\|_{L^\nf(B)}\le C\mko\eta,
\end{equation}
where $C=C(n,\al,p)$.
\end{Lm}
\begin{proof}
Set $G\coloneq g^{p-1}$. The definition~\eqref{eq:defLor} of Lorentz space implies
\begin{equation}\label{equivpowlor}
    \|G\|_{L^{\frac n{\al p},\frac{s}{p-1}}(B)}
    =\|g\|_{L^{q,s}(B)}^{p-1}\le\eta^{p-1}.
\end{equation}
We define $\mathbf W_{\alpha,p}G$ on $B$ by~\eqref{def:Wolffpot} and extend it by zero outside $B$. 
We first assume $s>1$. Combining~\eqref{comWolHaMaz} and~\eqref{equivpowlor} with Cianchi's estimate~\cite[Thm.~3.1(iii)]{Cianchi11} for the Havin--Maz'ya potential gives
\begin{equation}\label{eq:Cianest}
    \int_0^{|B|}
    \bigg(\frac{(\mathbf W_{\alpha,p}G)^*(\tau)}{1+\log(|B|/\tau)}\bigg)^s\,\frac{d\tau}{\tau}
    \le C\mko\eta^s,
\end{equation}
where $C=C(n,\al,p,s)$. Since $(\mathbf W_{\alpha,p}G)^*$ is non-increasing, it follows that for any $0<t\le |B|$,
\begin{align*}
    C\mko\eta^s
    &\ge \big(\mko(\mathbf W_{\alpha,p}G)^*(t)\mko\big)^s\int_0^t\big(1+\log(|B|/\tau)\big)^{-s}\,\frac{d\tau}{\tau}\\
    &=\frac1{s-1}\big(\mko(\mathbf W_{\alpha,p}G)^*(t)\mko\big)^s\big(1+\log(|B|/t)\big)^{1-s}.
\end{align*}
Hence, we have
\[
    (\mathbf W_{\alpha,p}G)^*(t)\le\wti C\mko\eta\mko\big(1+\log(|B|/t)\big)^{\frac{s-1}{s}},
\]
where $\wti C=\wti C(n,\al,p,s)$. Equimeasurability and Young's inequality then imply that
\begin{align*}
    \fint_B\exp\big(\ga \mathbf W_{\alpha,p}G\big)
    &\le\frac1{|B|}\int_0^{|B|}\exp\Big(\wti C\ga\eta\big(1+\log(|B|/t)\big)^{\frac{s-1}{s}}\Big)\,dt\\
    &=\int_0^\nf\exp\Big(\wti C\ga\eta(1+\si)^{\frac{s-1}{s}}-\si\Big)\,d\si\\
    &\le\int_0^\nf\exp\Big(C_4(\ga\eta)^s+\frac12-\frac\si2\Big)\,d\si\\
    &\le4\exp\big(C_4(\ga\eta)^s\big),
\end{align*}
where $C_4=C_4(n,\al,p,s)$.
This proves~\eqref{eq:WolLor} for $s>1$.

When $s=1$, the estimate~\eqref{eq:Wolend} follows from~\cite[Thm.~3.1(iv)]{Cianchi11} and~\eqref{equivpowlor}. 
Approximation gives $\mathbf W_{\alpha,p}G\in C^0(\ov B)$.
This completes the proof.
\end{proof}

\subsection{Interior reverse H\"older inequality}\label{sec:intgra}
\
\vskip5pt

The key step for our proof of the reverse H\"older inequality is a pointwise gradient estimate for homogeneous solutions obtained by adapting the proof of~\cite[Eq.~(1.9)]{Nguyen23} to our coefficient assumptions. 

In this subsection, we fix
\begin{equation}\label{setp1p0al}
   n\ge 2,\qquad  \frac{2n}{3}\le p_0<n,\qquad 1\le s<\nf.
\end{equation}
For a ball $B_{r}(x)\subset \Om$, we define
\begin{equation}\label{defvt_A}
     \vt_A(x,r):=r
     \bigg(\fint_{B_r(x)} |DA|^{p_0}\bigg)^{1/p_0}.
\end{equation}
We start with an $L^2$ gradient estimate for $L_A$.

\begin{Lm}\label{lm:L2grad}
 
Let $\Om\subset \R^n$ be a bounded Lipschitz domain, and let $B_{2r}(x)\subset \Om$. 
Suppose that the coefficient tensor $A=S+K$ satisfies~\eqref{bdsympart} and~\eqref{eq:intcond}, with $S$ and $K$ defined in~\eqref{def:spatparts}, and assume $\vt_A(x,2r)<\nf$. Let $f=(f^\al_i)\in L^{2}(B_{2r}(x),\mathbb R^{m\times n})$. 
Suppose $u\in W^{1,2}(B_{2r}(x),\mathbb R^m)$ satisfies
\[
\p_\al
\big(A^{\alpha\beta}_{ij}\partial_\beta u^j\big)
=\p_\al f^\al_i
\qquad\text{in }B_{2r}(x).
\]
Then there exist a constant $C=C(n,m,\La)>0$ such that
\begin{equation}\label{eq:L2grad}
\begin{aligned}
    \|Du\|_{L^2(B_r(x))}
    \le C\Big(\|f\|_{L^2(B_{2r}(x))}+ \big(1+\vt_A(x,2r)\mko\big)^{n+2}
    r^{-1-\frac n2}\|u\|_{L^1(B_{2r}(x))}
    \Big).
\end{aligned}
\end{equation}

\end{Lm}

\begin{proof}
By translation and rescaling, we may assume $x=0$ and $r=1$. Let $\{\bar K^{\al\beta}_{ij}\}$ be an arbitrary constant coefficient tensor satisfying $\bar K^{\al\beta}_{ij}=-\bar K^{\beta\al}_{ij}$. For any $\vp\in W^{1,2}_0(B_2,\R^m)$, integration by parts implies
\begin{equation*}
    \int_{B_2} \bar K^{\al\beta}_{ij}\mko \p_\beta u^j \p_\al \vp^i=\frac 12  \int_{B_2} \bar K^{\al\beta}_{ij}\mko \p_\beta u^j \p_\al \vp^i-\frac 12  \int_{B_2} \bar K^{\al\beta}_{ij}\mko\p_\al u^j \p_\beta \vp^i=0.
\end{equation*}
In other words, we have 
\begin{equation}\label{consantivan}
    \p_\al\big(\bar K^{\al\beta}_{ij}\p_\beta u^j\big)=0.
\end{equation}
Thus, we can assume without loss of generality that $(K^{\al\beta}_{ij})_{B_2}=0$ for any $i,j,\al,\beta$. 
Since $p_0\ge 2n/3$, Poincar\'e's inequality together with~\eqref{bdsympart} then implies 
    \begin{equation}\label{L2nnormA}
        \|A\|_{L^{2n}(B_2)}\le C(1+\vt_A(0,2)\mko\big),
    \end{equation}
where $C=C(\La,n)$.
Let $1\le\rho<\sigma\le2$. We choose
$\chi\in C_c^\infty(B_\sigma)$ such that $0\le\chi\le 1$ and
$|D\chi|\le C/(\sigma-\rho)$ on $B_\si$, and
$\chi=1$ on $B_\rho$.
Testing the equation with $\chi^2 u$ and applying H\"older's inequality yields
\begin{align}\label{intequa}
\begin{aligned}
   & \int_{B_\si} A^{\al\beta}_{ij} \p_\beta (\chi u^j) \p_\al(\chi u^i) \\
    &=\int_{B_\si} f_i^\al \p_\al (\chi^2 u^i) +\chi A^{\al\beta}_{ij} (\p_\al u^i \p_\beta \chi u^j-\p_\beta u^j \p_\al\chi\mko u^i) +A^{\al\beta}_{ij}\p_\beta \chi \p_\al \chi  \mko u^j u^i\\
    &\le  C\lf(\|f\|_{L^2(B_\si)}+(\si-\rho)^{-1}\|A\|_{L^{2n}(B_\si)}\|u\|_{L^{\frac{2n}{n-1}}(B_\si)}\rg)\lf(\|Du\|_{L^2(B_\si)}+\frac{ \|u\|_{L^2(B_\si)}}{\si-\rho}\rg).
    \end{aligned}
\end{align}
By the Gagliardo--Nirenberg inequality, we have
\begin{equation}\label{gagliNirenapp}
    \|u\|^2_{L^{\frac{2n}{n-1}}(B_\si)}\le C \mko\|u\|_{L^2(B_\si)} \big(\|Du\|_{L^{2}(B_\si)} +\|u\|_{L^{2}(B_\si)}\big).
\end{equation}
Let $\vae>0$ be arbitrary. Combining~\eqref{L2nnormA}--\eqref{gagliNirenapp} with Young's inequality, we obtain
\[
    \int_{B_\si} A^{\al\beta}_{ij} \p_\beta (\chi u^j) \p_\al(\chi u^i)\le \vae \|Du\|_{L^2(B_\si)}^2+\frac{C_\vae\big(1+\vt_A(0,2)\mko\big)^4}
{(\sigma-\rho)^4}\mko\|u\|_{L^2(B_\si)}^2  +C_\vae \|f\|_{L^{2}(B_2)}^2,
\]
where $C_\vae$ depends only on $\vae,\La, n, m$. 
Set $\vae=(4\La)^{-1}$. By~\eqref{intequa} and the integral ellipticity condition~\eqref{eq:intcond}, we then obtain
\begin{equation}\label{firiteine}
    \|Du\|_{L^2(B_\rho)}^2\le \frac 14\mko\|Du\|_{L^2(B_\si)}^2+\frac{C\big(1+\vt_A(0,2)\mko\big)^4}
{(\sigma-\rho)^4}\mko \|u\|_{L^2(B_\si)}^2
+C
\|f\|_{L^{2}(B_2)}^2,
\end{equation}
where $C=C(\La,n,m)$.
Applying the Gagliardo--Nirenberg inequality again implies
\[
\|u\|_{L^2(B_\sigma)}
\le C\mko \|u\|_{L^1(B_\si)}^{\frac 2{n+2}} \big(\|Du\|_{L^{2}(B_\si)}+\|u\|_{L^1(B_\si)}\big)^{\frac{n}{n+2}}\le C\|u\|_{L^1(B_\si)}^{\frac 2{n+2}}  \|Du\|_{L^2(B_\si)}^{\frac{n}{n+2}}+C\|u\|_{L^1(B_\si)}.
\]
Combining this with~\eqref{firiteine}
and applying Young's inequality, we obtain that
\[
    \|Du\|_{L^2(B_\rho)}^2
    \le \frac 12 \|Du\|_{L^2(B_\si)}^2
    +\frac{C\big(1+\vt_A(0,2)\mko\big)^{2n+4}}
    {(\sigma-\rho)^{2n+4}}\mko\|u\|_{L^1(B_\si)}^2+C\|f\|_{L^{2}(B_2)}^2,
\]
where $C=C(\La,n,m)$.
The estimate~\eqref{eq:L2grad} then follows from standard iteration; see~\cite[Lem.~1.1]{GiaGiu82}.
\end{proof}

Following the proof of~\cite[Lem.~2.1]{Nguyen23}, we obtain the following harmonic replacement estimate. 
The argument requires distinguishing large and small values of $\vt_A(x,2r)$.
\begin{Lm}\label{lm:intcom}
Let $\Om$, $B_{2r}(x)$, and $A$ be as in Lemma~\ref{lm:L2grad}. Set $\bar A\coloneq (A)_{B_{2r}(x)}$. Suppose $w,h\in W^{1,2}(B_{2r}(x),\R^m)$ satisfy
\begin{equation}\label{eqsuf}
    \begin{dcases}
    \begin{aligned}
        L_A w&\,=\, 0&& \text{ in }\;B_{2r}(x),\\
        L_{\bar A}h&\,=\,0&& \text{ in }\; B_{2r}(x),\\
        h&\,=\,w&& \text{ on }\; \p B_{2r}(x).
        \end{aligned}
    \end{dcases}
\end{equation}
Then there exists a constant $C=C(n,m,\La)$ such that
\begin{equation}\label{eq:intcom}
 \|D(w-h)\|_{L^2(B_{3r/2}(x))}\le C \mko  \vt_A(x,2r)\|Dw\|_{L^2(B_{2r}(x))}.
\end{equation}
\end{Lm}
\begin{proof}
    By~\eqref{consantivan}, we can assume without loss of generality that $(K)_{B_{2r}(x)}=0$, and hence $\bar A=\bar S$, where $\bar S=(S)_{B_{2r}(x)}$.
    By the $L^\nf$ bound~\eqref{bdsympart} and the Legendre--Hadamard condition~\eqref{eq:LegHadcon}, the standard elliptic estimates apply to the constant-coefficient operator $L_{\bar A}=L_{\bar S}$; see for instance~\cite{DouNir55,AgmDouNir64,GiaqMari12}.

    By~\eqref{eqsuf}, we have
    \begin{equation}\label{eqw-hLAw}
    \begin{dcases}
    \begin{aligned}
        L_{\bar A}(w-h)&\,=\,L_{\bar A} w&& \text{ in }\; B_{2r}(x),\\
        w-h&\,=\,0&& \text{ on }\; \p B_{2r}(x).
        \end{aligned}
    \end{dcases}
\end{equation}
The $L^2$ estimate then gives
\begin{equation}\label{eqintcomrou}
     \|D(w-h)\|_{L^2(B_{2r}(x))}\le C  \|L_{\bar A} w\|_{W^{-1,2}(B_{2r}(x))}\le C\|Dw\|_{L^2(B_{2r}(x))}.
\end{equation}
This proves~\eqref{eq:intcom} in the case $\vt_A(x,2r)\ge 1$. 

Now we assume $\vt_A(x,2r)< 1$. 
By~\eqref{eqintcomrou} and using the equation $L_{\bar A} h=0$, we have
\begin{equation}\label{DhLnfbd}
    \|Dh\|_{L^\nf(B_{7r/4}(x))}\le C\mko r^{-\frac n2}\|Dh\|_{L^2(B_{2r}(x))}\le C\mko r^{-\frac n2}\|Dw\|_{L^2(B_{2r}(x))}.
\end{equation}
We write $h=(h^j)_{1\le j\le m}$. The systems in~\eqref{eqsuf} give
\[
    \p_\al
\Big(A^{\alpha\beta}_{ij}\partial_\beta \big(w^j-h^j\big)\Big)=\p_\al
\Big(\big(\bar A^{\al\beta}_{ij}-A^{\alpha\beta}_{ij}\big)\partial_\beta h^j\Big).
\]
Applying Lemma~\ref{lm:L2grad} together with a covering and rescaling argument then yields 
\begin{align}\label{w-hL2finest}
\begin{aligned}
     &\|D(w-h)\|_{L^2(B_{3r/2}(x))}\\
     &\le C\Big( \|\bar A-A\|_{L^2(B_{2r}(x))}\|Dh\|_{L^\nf (B_{7r/4}(x))}+r^{-1-\frac n2}\|w-h\|_{L^1(B_{2r}(x))}\Big),
     \end{aligned}
\end{align}
where $C=C(n,m,\La)$. Poincar\'e's inequality gives
\begin{equation}\label{Poinoscder}
    r^{-\frac n2}\|\bar A-A\|_{L^2(B_{2r}(x))}\le C\mko \vt_A(x,2r),
\end{equation}
where $C$ depends only on $n,m,\La$ and is independent of $p_0$.
To estimate $\|w-h\|_{L^1(B_{2r}(x))}$, we use the following system
\begin{equation*}
    \begin{dcases}
    \begin{aligned}
        \p_\al
\Big(\bar A^{\alpha\beta}_{ij}\partial_\beta \big(w^j-h^j\big)\Big)&\,=\,\p_\al
\Big(\big(\bar A^{\al\beta}_{ij}-A^{\alpha\beta}_{ij}\big)\partial_\beta w^j\Big)&& \text{ in }\; B_{2r}(x),\\
        w-h&\,=\,0&& \text{ on }\; \p B_{2r}(x).
        \end{aligned}
    \end{dcases}
\end{equation*}
The standard $W^{2,p}$ estimate (see~\cite{AgmDouNir64}), combined with Morrey's embedding, implies that for any $f=(f^i)\in L^\nf (B_{2r}(x),\R^m)$, the dual problem 
\begin{equation*}
    \begin{dcases}
    \begin{aligned}
        \p_\al
    \Big(\bar A^{\beta\al}_{ji}\partial_\beta \phi^j\Big)&\,=\,f^i&& \text{ in }\; B_{2r}(x),\\
        \phi&\,=\,0&& \text{ on }\; \p B_{2r}(x),
        \end{aligned}
    \end{dcases}
\end{equation*}
admits a unique solution $\phi=(\phi^j)\in W^{1,\nf}(B_{2r}(x),\R^m)$ with
\[
    \|D\phi\|_{L^\nf(B_{2r}(x))}\le C\mko r\|f\|_{L^\nf(B_{2r}(x))}.
\]
By duality, we obtain
\begin{equation}\label{w-hL1est}
    \|w-h\|_{L^1(B_{2r}(x))}\le C\mko r\mko\|\bar A-A\|_{L^2(B_{2r}(x))}\|Dw\|_{L^2(B_{2r}(x))}.
\end{equation}
The estimate~\eqref{eq:intcom} follows by combining~\eqref{w-hL2finest}--\eqref{w-hL1est}.
\end{proof}
Applying Lemma~\ref{lm:intcom} and adapting the proof of~\cite[Eq.~(1.9)]{Nguyen23} leads to the following gradient estimate. 
Here we control the coefficient term pointwise by a Wolff potential without imposing a Dini mean oscillation condition.
We note that the argument in~\cite{Nguyen23} has an oversight error at equation (2.12) therein, which can be corrected by a minor modification that does not affect the proof. 
See~\eqref{hk-hk+1preest} below.
\begin{Lm}\label{lm:Nguint}
Let $\Om$ and $A$ be as in Lemma~\ref{lm:L2grad}. Let $B_{R}(x)\subset \Om$, $0<r\le R/4$, and we assume $\vt_A(x,R)<\nf$. Let $w\in W^{1,2}(B_{R}(x),\R^m)$ solve $L_A w=0$.
Then there exists a constant $C=C(n,m,\La)>0$ such that
\begin{align}\label{Nguyenptbd}
 \|Dw\|_{L^2(B_r(x))}
 &\le C\Big(\frac{r}R\Big)^{n/2}\exp\bigg(C\int_r^R \vt_A(x,t)\,\frac{dt}{t}\bigg)
  \|Dw\|_{L^2(B_R(x))}.
\end{align}
\end{Lm}
\begin{proof}
By translation, we may assume $x=0$. For $k\ge0$, let $r_k=4^{-k} R/2$, $\bar A_k=(A)_{B_{2r_k}}$, and $h_k$ be the solution of
  \begin{equation}\label{eq:LHhkw}
    \begin{dcases}
    \begin{aligned}
        L_{\bar A_k}h_k&\,=\,0&& \text{ in }\; B_{2r_k},\\
        h_k&\,=\,w&& \text{ on }\; \p B_{2r_k}.
        \end{aligned}
    \end{dcases}
\end{equation}
We set
\[
    a_k=r_k^{-n/2}\|D(w-h_k)\|_{L^2(B_{r_k})},\qquad b_k=\|Dh_k\|_{L^\nf(B_{r_k})}.
\]
The triangle inequality gives
\begin{align}\label{trianinedw}
    \|Dw\|_{L^2(B_{r_k})}\le C \mko r_k^{n/2}(a_k+b_k),
\end{align}
where $C=C(n)$.
By Lemma~\ref{lm:intcom}, we then obtain
\begin{align}\label{estak+1}
\begin{aligned}
    a_{k+1}&\le r_{k+1}^{-n/2}\|D(w-h_{k+1})\|_{L^2(B_{3r_{k+1}/2})}\\
    &\le C\mko  r_{k+1}^{-n/2}\mko \vt_A(x,2r_{k+1})\|Dw\|_{L^2(B_{2r_{k+1}})}\\
    &\le C\mko \vt_A(x,2r_{k+1}) (a_k+b_k),
\end{aligned}
\end{align}
where $C=C(n,m,\La)$. Let $k\ge 1$.
It follows from~\eqref{estak+1} that 
\begin{align}\label{hk-hk+1preest}
\begin{aligned}
    r_k^{-n/2}\|D(h_k-h_{k+1})\|_{L^2(B_{3r_{k+1}/2})}&\le a_k+r_k^{-n/2}\|D(w-h_{k+1})\|_{L^2(B_{3r_{k+1}/2})}\\
    &\le C\mko \vt_A(x,2r_{k})\mko(a_{k-1}+b_{k-1}+a_k+b_k).
    \end{aligned}
\end{align}
By~\eqref{eq:LHhkw}, we have
\begin{equation}\label{eqhk-hk+1}
    L_{\bar A_k} (h_k-h_{k+1})=(L_{\bar A_{k+1}}-L_{\bar A_k})\mko h_{k+1}.
\end{equation}
As in the proof of Lemma~\ref{lm:intcom}, standard elliptic estimates for~\eqref{eq:LHhkw} imply that
\begin{align*}
    &r_{k+1}^{n/2}\Big(\|Dh_{k+1}\|_{L^\nf(B_{3r_{k+1}/2})}+ r_{k+1}\|D^2h_{k+1}\|_{L^\nf(B_{3r_{k+1}/2})}\Big)\\[1ex]
   &\le C\|Dh_{k+1}\|_{L^2(B_{2r_{k+1}})}\le C\|Dw\|_{L^2(B_{2r_{k+1}})}.
\end{align*}
Hence, by Poincar\'e's inequality and~\eqref{trianinedw}, we obtain
\begin{align}\begin{aligned}\label{hk-hk+1Lnfbd}
    \big\|(L_{\bar A_{k+1}}-L_{\bar A_k})\mko h_{k+1}\big\|_{L^\nf(B_{3r_{k+1}/2})}&\le C \mko|\bar A_{k+1}-\bar A_k|\mko\|D^2h_{k+1}\|_{L^\nf(B_{3r_{k+1}/2})}\\
    &\le C\mko \vt_A(x,2r_{k})\mko r_k^{-n/2-1} \|Dw\|_{L^2(B_{2r_{k+1}})}\\
    &\le C\mko \vt_A(x,2r_{k})\mko r_k^{-1}(a_k+b_k).
    \end{aligned}
\end{align}
Applying the interior $W^{2,p}$ estimate (see e.g.~\cite{DouNir55}) to~\eqref{eqhk-hk+1} and using Morrey's embedding, we deduce
\begin{align*}
    &\|D(h_k-h_{k+1})\|_{L^\nf(B_{r_{k+1}})}\\
    &\le C \Big(r_{k+1}\mko  \big\|(L_{\bar A_{k+1}}-L_{\bar A_k})\mko h_{k+1}\big\|_{L^\nf(B_{3r_{k+1}/2})}+   r_{k+1}^{-n/2}\|D(h_k-h_{k+1})\|_{L^2(B_{3r_{k+1}/2})}\Big)\\
    &\le  C\mko \vt_A(x,2r_{k})\mko(a_{k-1}+b_{k-1}+a_k+b_k).
\end{align*}
Combined with~\eqref{estak+1} and triangle inequality, this implies
\begin{equation}\label{ak+bkiteine}
    a_{k+1}+b_{k+1}\le b_k+C\mko \vt_A(x,2r_{k})\mko(a_{k-1}+b_{k-1}+a_k+b_k).
\end{equation}
Set 
\[
    M_k=\max_{0\le \ell\le k}(a_\ell+b_\ell).
\]
The inequality~\eqref{ak+bkiteine} then yields
\[
    M_{k+1}\le (1+C\mko \vt_A(x,2r_{k})) M_k,\qquad \forall\, k\ge 1.
\]
In addition, the proof of~\eqref{DhLnfbd} gives 
\[
    M_1\le C\mko R^{-n/2}\|Dw\|_{L^2(B_{R})},
\]
where $C=C(n,m,\La)$. It follows that
\begin{align*}
     r_{k+1}^{-n/2}\|Dw\|_{L^2(B_{r_{k+1}})}& \le C\mko M_{k+1}\le C \exp \lf(C\sum_{\ell=1}^{k}\vt_A(x,2r_{\ell})\rg)M_1\\
     &\le C\exp\bigg(C\int_{2r_k}^R \vt_A(x,t)\,\frac{dt}{t}\bigg) R^{-n/2} \|Dw\|_{L^2(B_R)},
\end{align*}
where $C=C(n,m,\La)$.
This completes the proof.
\end{proof}
We now obtain the interior reverse H\"older inequality.
\begin{Co}\label{co:intrevHol}
Suppose $\Om$ and $A$ satisfy the conditions of Theorem~\ref{th:divest}. Let $B_{2R}(x_0)\subset \Om$, $2<q<\nf$, and let $w\in W^{1,2}(B_{2R}(x_0),\R^m)$ solve $L_A w=0$. Then there exist positive constants $C_5=C_5(n,m,\La,s)$ and $C=C(n,m,\La)$ such that
\begin{equation}\label{eq:intrevHol}
 \bigg(\fint_{B_{R}(x_0)}|Dw|^q\bigg)^{1/q}
 \le C\exp(C_5 \mko q^{s-1}\mko \eta^s)\bigg(\fint_{B_{2R}(x_0)}|Dw|^2\bigg)^{1/2}.
\end{equation}
\end{Co}
\begin{proof}
Set $p_0=2n/3$, $G=|DA|^{p_0}\mkt\mathbf 1_{B_{2R}(x_0)}$, $p_1\coloneq p_0+1$, and $\alpha_0\coloneq p_0/p_1$. 
Letting $r\to 0$ in~\eqref{Nguyenptbd} and applying the Lebesgue differentiation theorem, we obtain that for almost every $x\in B_R(x_0)$,
\begin{align*}
    |Dw(x)|&\le C\exp\lf(C\int_0^R \vt_A(x,t)\,\frac{dt}{t}\rg)
 \bigg(\fint_{B_R(x)}|Dw|^2\bigg)^{1/2}\\
    & \le C\exp \Big(C\mko W_{\alpha_0,p_1}G(x)\mko\Big)\bigg(\fint_{B_{2R}(x_0)}|Dw|^2\bigg)^{1/2},
\end{align*}
where $C=C(n,m,\La)$.
It follows that 
\[
     \bigg(\fint_{B_{R}(x_0)}|Dw|^q\bigg)^{1/q}
 \le C \bigg(\fint_{B_R(x_0)} \exp \Big(C\mko q\mko W_{\alpha_0,p_1}G\mko\Big)\bigg)^{\frac 1q}\bigg(\fint_{B_{2R}(x_0)}|Dw|^2\bigg)^{1/2}.
\]
Applying Lemma~\ref{lm:Wolff-exp} with $g=|DA|\mko \mathbf1_{B_{2R}(x_0)}$ completes the proof of~\eqref{eq:intrevHol}.
\end{proof}

\section{Global estimates}\label{sec:global}
In this section, we prove the boundary reverse H\"older inequality and complete the proofs of Theorems~\ref{th:divest} and~\ref{Th:divHajespa}. Set $\R^n_+\coloneq \{(x',x_n)\in \R^n: x_n> 0\}$.
For $x\in \p \R^n_+$, we denote 
\[
   B_r^+(x)\coloneq B_r(x)\cap \R_+^n,\qquad \Ga_r(x)\coloneq B_r(x)\cap \p \R^n_+,\qquad B_r^+\coloneq B_r^+(0),\qquad \Ga_r\coloneq \Ga_r(0).
\]
Throughout this section, let $\Om$ be bounded and $(\de,R_0)$-Lipschitz for some $R_0>0$, where $\de\in (0,\frac 1{10})$ is to be determined. 
For $x\in \ov\Om$ and $r>0$, we write $\Om_r(x)\coloneq B_r(x)\cap \Om$.
We fix $x_0\in \p\Om$.  
By Definition~\ref{def:smaLip}, there exists a bi-Lipschitz map $\Psi$ from $\ov{\Om_{R_0}(x_0)}$ to $\ov{\R^n_+}$ such that $\Psi(x_0)=0$, $\det(D\Psi)=1$, and
\begin{equation}\label{propbiLip}
 \Psi\big(\p\Om\cap B_{R_0}(x_0) \big)=\Psi\big(\ov{\Om_{R_0}(x_0)}\big)\cap \p \R^n_+,\qquad\|D\Psi-I_n\|_{L^\nf(\Om_{R_0}(x_0))}\le\delta.
\end{equation}
 In particular, this implies that for any $|x-x_0|=r$ with $0<r\le R_0$, we have $3r/4\le |\Psi(x)| \le  5r/4$. By shrinking $R_0$ if necessary, we can assume without loss of generality that
 \begin{equation}\label{incimage}
    B_{3r/4}^+ \subset  \Psi (\Om_r(x_0))\subset B_{5r/4}^+,\qquad \forall\, 0<r\le R_0.
 \end{equation}
Until the proof of Theorem~\ref{Th:divHajespa}, we assume the coefficient tensor $A$ satisfies the conditions of Theorem~\ref{th:divest}. On $B^+_{3R_0/4}$, we define $\Phi\coloneq \Psi^{-1}$ and
\begin{equation}\label{deftiA}
   \wti A^{\al\beta}_{ij}\coloneq \big(\mko(D\Psi)_{\al\ga} A^{\ga\si}_{ij} (D\Psi)_{\beta\si}\big)\circ \Phi,
\end{equation}
where we write $\Psi=(\Psi_\al)_{1\le \al\le n}$ and $(D\Psi)_{\al\ga}\coloneq \p_\ga \Psi_\al$.
Then for any $0<r\le R_0$ and $w\in W^{1,2}(\Om_{r}(x_0),\R^m)$ solving $L_Aw=0$, the function $\wti w =w \circ \Phi$ satisfies
\begin{equation}\label{transsat}
    (L_{\wti A}\wti w )^i=\p_\al\big( \wti A^{\al\beta}_{ij}\p_\beta \wti w^j\big)=0 \quad \text{ in }B_{3r/4}^+.
\end{equation}
As in Section~\ref{sec:intgra}, we set $2n/3\le p_0<n$, and for all $x\in \ov\Om\cap B_{R_0}(x_0)$, $0<r\le R_0$, we define \begin{equation}\label{defvt_Abdy}
     \vt_A(x,r):=r
     \bigg(\fint_{\Om_r(x)} |DA|^{p_0}\bigg)^{1/p_0}.
\end{equation}
\subsection{Boundary reverse H\"older inequality}
\
\vskip5pt
We prove an analog of~\eqref{eq:L2grad} for boundary estimates.
\begin{Lm}\label{L2firgrabdy}
    Let $z\in \Ga_{R_0/4}$, $r<R_0/8$, and let $f=(f^\al_i)\in L^{2}(B^+_{2r}(z),\R^{m\times n})$.
Suppose $u\in W^{1,2}(B_{2r}^+(z),\R^m)$ with $u=0$ on $\Ga_{2r}(z)$ satisfies
\[
\p_\al
\big(\wti A^{\alpha\beta}_{ij}\partial_\beta u^j\big)
=\p_\al f^\al_i
\qquad\text{in }B_{2r}^+(z).
\]
Then there exists a constant $C=C(n,m,\La)>0$ such that
\begin{equation}\label{eq:L2gradbry}
\begin{aligned}
    \|Du\|_{L^2(B_r^+(z))}
    \le C\Big(\|f\|_{L^2(B_{2r}^+(z))}
    +\big(1+\vt_A(\Phi(z),3r)\big)^{n+2}
    r^{-1-\frac n2}\|u\|_{L^1(B_{2r}^+(z))}\Big).
\end{aligned}
\end{equation}
\end{Lm}

\begin{proof}
By~\eqref{consantivan} and~\eqref{transsat}, we can assume without loss of generality that
\[
    (K\circ\Phi)_{B_{2r}^+(z)}=0.
\]
Poincar\'e's inequality together with~\eqref{bdsympart} then implies
\begin{equation}\label{eq:bdycoeffL2n}
r^{-\frac12}\|\wti A\|_{L^{2n}(B_{2r}^+(z))}
\le C\big(1+\vt_A(\Phi(z),3r)\big),
\end{equation}
where $C=C(n,m,\La)$.
Let $r\le\rho<\si\le2r$. We choose $\chi\in C_c^\nf(B_\si(z))$ such that $0\le\chi\le 1$ and $|D\chi|\le C/(\si-\rho)$ on $B_\si (z)$, and $\chi=1$ on $B_\rho(z)$. Since $u=0$ on $\Ga_{2r}(z)$, we can test the equation with $\chi^2u$. The proof of Lemma~\ref{lm:L2grad}, using~\eqref{eq:bdycoeffL2n} in place of~\eqref{L2nnormA}, gives
\begin{align*}
 &\|Du\|_{L^2(B_\rho^+(z))}^2\\
 &\le \frac12\|Du\|_{L^2(B_\si^+(z))}^2
 +\frac{Cr^{n+2}\big(1+\vt_A(\Phi(z),3r)\big)^{2n+4}}{(\si-\rho)^{2n+4}}\|u\|_{L^1(B_\si^+(z))}^2+C\|f\|_{L^2(B_{2r}^+(z))}^2,
\end{align*}
where $C=C(n,m,\La)$. The estimate~\eqref{eq:L2gradbry} follows from~\cite[Lem.~1.1]{GiaGiu82}.
\end{proof}

To obtain the boundary harmonic replacement estimate, we fix a smooth domain $D\subset\R^n_+$ such that
\begin{equation}\label{eq:capin}
 B_{7/4}^+\subset D\subset B_2^+\qquad\text{and}\qquad \p D\cap\p\R^n_+\supset\Ga_{7/4}.
\end{equation}
For $z\in\p\R^n_+$ and $r>0$, set $D_r(z)\coloneq z+rD$.
Following the proof of Lemma~\ref{lm:intcom}, we obtain the following boundary harmonic replacement estimate.
\begin{Lm}\label{lm:bdycom}
Let $z\in\Ga_{R_0/4}$ and $0<r<R_0/8$. Set $\bar A\coloneq(A\circ\Phi)_{D_r(z)}$. Suppose $\wti w,h\in W^{1,2}(D_r(z),\R^m)$ satisfy
\begin{equation}\label{eq:bdysys}
 \begin{dcases}
 \begin{aligned}
 L_{\wti A}\wti w&\,=\,0&&\text{ in }\;D_r(z),\\
 L_{\bar A}h&\,=\,0&&\text{ in }\;D_r(z),\\
 h&\,=\,\wti w&&\text{ on }\;\p D_r(z),
 \end{aligned}
 \end{dcases}
\end{equation}
and assume $\wti w=0$ on $\Ga_{7r/4}(z)$. Then there exists a constant $C=C(n,m,\La)>0$ such that
\begin{equation}\label{eq:bdycom}
 \|D(\wti w-h)\|_{L^2(B_{3r/2}^+(z))}
 \le C\big(\vt_A(\Phi(z),3r)+\de\big)\|D\wti w\|_{L^2(D_r(z))}.
\end{equation}
\end{Lm}
\begin{proof}
 By~\eqref{consantivan}, we have $L_{\bar A}=L_{\bar S}$, where $\bar S=(S\circ\Phi)_{D_r(z)}$. As in Lemma~\ref{L2firgrabdy}, we can assume without loss of generality that $(K\circ\Phi)_{D_r(z)}=0$.
 Combining~\eqref{bdsympart} and~\eqref{propbiLip} with Poincar\'e's inequality gives
\begin{equation}\label{eq:comcoeffbd}
 \|\wti A-\bar S\|_{L^2(D_r(z))}
 \le Cr^{n/2}\big(\vt_A(\Phi(z),3r)+\de\big),
\end{equation}
where $C=C(n,m,\La)$. Since $L_{\bar S}(\wti w-h)=L_{\bar S} \wti w$ in $D_{r}(z)$ and $h=\wti w$ on $\p D_{r}(z)$, the $L^2$ estimate gives
\[
 \|D(\wti w-h)\|_{L^2(D_r(z))}\le C\|D\wti w\|_{L^2(D_r(z))}.
\]
Since $B^+_{3r/2}(z)\subset D_r(z)$, this proves~\eqref{eq:bdycom} when $\vt_A(\Phi(z),3r)\ge1$.

Now we assume $\vt_A(\Phi(z),3r)<1$.
Standard elliptic estimates in~\cite{AgmDouNir64} imply that
\begin{equation}\label{eq:DhLinfbdy}
 \|Dh\|_{L^\nf(B_{13r/8}^+(z))}
 \le Cr^{-n/2}\|Dh\|_{L^2(D_r(z))}
 \le Cr^{-n/2}\|D\wti w\|_{L^2(D_r(z))}.
\end{equation}
The system~\eqref{eq:bdysys} gives
\[
 \p_\al\Big(\wti A_{ij}^{\al\beta}\p_\beta(\wti w^j-h^j)\Big)
 =\p_\al\Big(\big(\bar S_{ij}^{\al\beta}-\wti A_{ij}^{\al\beta}\big)\p_\beta h^j\Big).
\]
Since $B_{3r/2}^+(z)\subset B_{13r/8}^+(z)\subset D_r(z)$, applying Lemma~\ref{L2firgrabdy} together with a covering and rescaling argument yields
\begin{align}\label{eq:bdycomL2}
\begin{aligned}
 &\|D(\wti w-h)\|_{L^2(B_{3r/2}^+(z))}
 \\&\le C\Big(\|\wti A-\bar S\|_{L^2(D_r(z))}\|Dh\|_{L^\nf(B_{13r/8}^+(z))}+r^{-1-n/2}\|\wti w-h\|_{L^1(D_r(z))}\Big).
\end{aligned}
\end{align}
To estimate the last term, we use
\[
 \begin{dcases}
 \begin{aligned}
 L_{\bar S}(\wti w-h)&\,=\,\p_\al\Big(\big(\bar S_{ij}^{\al\beta}-\wti A_{ij}^{\al\beta}\big)\p_\beta \wti w^j\Big)&&\text{ in }\;D_r(z),\\
 \wti w-h&\,=\,0&&\text{ on }\;\p D_r(z).
 \end{aligned}
 \end{dcases}
\]
For $f=(f^i)\in L^\nf(D_r(z),\R^m)$, the standard $W^{2,p}$ estimate in~\cite{AgmDouNir64}, combined with Morrey's embedding, gives a unique solution $\phi$ to the dual problem
\begin{equation*}
    \begin{dcases}
    \begin{aligned}
        \p_\al
    \Big(\bar A^{\beta\al}_{ji}\partial_\beta \phi^j\Big)&\,=\,f^i&& \text{ in }\; D_{r}(z),\\
        \phi&\,=\,0&& \text{ on }\; \p D_{r}(z),
        \end{aligned}
    \end{dcases}
\end{equation*}
with
\[
 \|D\phi\|_{L^\nf(D_r(z))}\le Cr\|f\|_{L^\nf(D_r(z))}.
\]
By duality, we obtain
\begin{equation}\label{eq:bdycomL1}
 \|\wti w-h\|_{L^1(D_r(z))}
 \le Cr\|\bar S-\wti A\|_{L^2(D_r(z))}\|D\wti w\|_{L^2(D_r(z))}.
\end{equation}
Combining~\eqref{eq:comcoeffbd}--\eqref{eq:bdycomL1} proves~\eqref{eq:bdycom}.
\end{proof}

Applying Lemma~\ref{lm:bdycom} and adapting the proof of Lemma~\ref{lm:Nguint}, we obtain the following boundary gradient estimate.
\begin{Lm}\label{lm:Ngubdy}
Let $z\in\Ga_{R_0/4}$, $0<R<R_0/8$, and $0<r\le R/4$. Suppose $\wti w\in W^{1,2}(B_R^+(z),\R^m)$ solves $L_{\wti A}\wti w=0$ with $\wti w=0$ on $\Ga_R(z)$. Then there exists a constant $C=C(n,m,\La)>0$ such that
\begin{equation}\label{eq:Ngubdy}
 \|D\wti w\|_{L^2(B_r^+(z))}
 \le C\Big(\frac rR\Big)^{n/2-C\de}
 \exp\bigg(C\int_r^R\vt_A(\Phi(z),t)\,\frac{dt}{t}\bigg)
 \|D\wti w\|_{L^2(B_R^+(z))}.
\end{equation}
\end{Lm}

\begin{proof}
For $k\ge0$, let $r_k=4^{-k-1}R$, $\bar A_k=(A\circ\Phi)_{D_{r_k}(z)}$, and $h_k$ be the solution of
\begin{equation}\label{eq:bdyrpl}
 \begin{dcases}
 \begin{aligned}
 L_{\bar A_k}h_k&\,=\,0&&\text{ in }\;D_{r_k}(z),\\
 h_k&\,=\,\wti w&&\text{ on }\;\p D_{r_k}(z).
 \end{aligned}
 \end{dcases}
\end{equation}
We set
\[
 a_k=r_k^{-n/2}\|D(\wti w-h_k)\|_{L^2(B_{r_k}^+(z))},\qquad
 b_k=\|Dh_k\|_{L^\nf(B_{r_k}^+(z))}.
\]
The triangle inequality gives
\begin{equation}\label{eq:bdyrtri}
 \|D\wti w\|_{L^2(B_{r_k}^+(z))}\le Cr_k^{n/2}(a_k+b_k).
\end{equation}
Since $B_{r_{k+1}}^+(z)\subset B_{3r_{k+1}/2}^+(z)$ and $D_{r_{k+1}}(z)\subset B_{r_k}^+(z)$, Lemma~\ref{lm:bdycom} implies
\begin{align}\label{eq:bdyak}
\begin{aligned}
 a_{k+1}
 &\le r_{k+1}^{-n/2}\|D(\wti w-h_{k+1})\|_{L^2(B_{3r_{k+1}/2}^+(z))}\\
 &\le Cr_{k+1}^{-n/2}\big(\vt_A(\Phi(z),3r_{k+1})+\de\big)\|D\wti w\|_{L^2(D_{r_{k+1}}(z))}\\
 &\le C\big(\vt_A(\Phi(z),3r_{k+1})+\de\big)(a_k+b_k),
\end{aligned}
\end{align}
where $C=C(n,m,\La)$.
Let $k\ge1$. Using~\eqref{eq:bdyak} and $\vt_A(\Phi(z),3r_{k+1})\le C\vt_A(\Phi(z),3r_k)$, we obtain
\begin{align*}
 r_{k}^{-n/2}\|D(h_k-h_{k+1})\|_{L^2(B_{3r_{k+1}/2}^+(z))}&\le a_k+r_k^{-n/2}\|D(\wti w-h_{k+1})\|_{L^2(B_{3r_{k+1}/2}^+(z))}\\
 &\le C\big(\vt_A(\Phi(z),3r_k)+\de\big)
 (a_{k-1}+b_{k-1}+a_k+b_k).
\end{align*}
Moreover, writing $\bar S_k=(S\circ\Phi)_{D_{r_k}(z)}$, we have
\begin{equation}\label{eq:bdyhkdiff}
 L_{\bar S_k}(h_k-h_{k+1})
 =(L_{\bar S_{k+1}}-L_{\bar S_k})h_{k+1}.
\end{equation}
By Poincar\'e's inequality, $|\bar S_k-\bar S_{k+1}|\le C\vt_A(\Phi(z),3r_k)$. Since $L_{\bar A_k}=L_{\bar S_k}$, standard boundary estimates for~\eqref{eq:bdyrpl} imply that
\begin{align*}
    &r_{k+1}^{n/2}\Big(\|Dh_{k+1}\|_{L^\nf(B_{3r_{k+1}/2}^+(z))}+ r_{k+1}\|D^2h_{k+1}\|_{L^\nf(B_{3r_{k+1}/2}^+(z))}\Big)\\[1ex]
   &\le C\|Dh_{k+1}\|_{L^2(D_{r_{k+1}}(z))}\le C\|D\wti w\|_{L^2(D_{r_{k+1}}(z))}.
\end{align*}
Applying the boundary $W^{2,p}$ estimate (see e.g.~\cite{AgmDouNir64}) to~\eqref{eq:bdyhkdiff} and using Morrey's embedding, we deduce
\begin{align*}
 &\|D(h_k-h_{k+1})\|_{L^\nf(B_{r_{k+1}}^+(z))}\\
 &\le C \Big(r_{k+1}\mko  \big\|(L_{\bar S_{k+1}}-L_{\bar S_k})\mko h_{k+1}\big\|_{L^\nf(B_{3r_{k+1}/2}^+(z))}+   r_{k+1}^{-n/2}\|D(h_k-h_{k+1})\|_{L^2(B_{3r_{k+1}/2}^+(z))}\Big)\\
 &\le C\big(\vt_A(\Phi(z),3r_k)+\de\big)(a_{k-1}+b_{k-1}+a_k+b_k).
\end{align*}
Combined with~\eqref{eq:bdyak}, this gives
\begin{equation}\label{eq:bdyMrec}
 a_{k+1}+b_{k+1}\le b_k+C\big(\vt_A(\Phi(z),3r_k)+\de\big)
 (a_{k-1}+b_{k-1}+a_k+b_k).
\end{equation}
Set $M_k=\max_{0\le\ell\le k}(a_\ell+b_\ell)$. The estimate~\eqref{eq:bdyMrec} implies that for any $k\ge1$,
\[
 M_{k+1}\le\big(1+C\vt_A(\Phi(z),3r_k)+C\de\big)M_k,
\]
where $C=C(n,m,\La)$.
In addition, the proof of~\eqref{eq:DhLinfbdy} gives
\[
 M_1\le CR^{-n/2}\|D\wti w\|_{L^2(B_R^+(z))}.
\]
Therefore, for $k\ge1$, we obtain
\begin{equation}\label{eq:bdygrow}
 \begin{aligned}
 &r_k^{-n/2}\|D\wti w\|_{L^2(B_{r_k}^+(z))}\\
 &\le CM_k\le C\exp\bigg(C\sum_{\ell=1}^{k-1}\vt_A(\Phi(z),3r_\ell)+Ck\de\bigg) M_1
\\
 &\le C\exp\bigg(C\int_{3r_k}^{3r_0} \vt_A(\Phi(z),t)\,\frac{dt}{t}+C\mko \de \log \Big( \frac {R}{r_k}\Big)\bigg)  R^{-n/2}\|D\wti w\|_{L^2(B_R^+(z))}.
 \end{aligned}
\end{equation}
The estimate~\eqref{eq:Ngubdy} follows.
\end{proof}

We now obtain the boundary reverse H\"older inequality.
\begin{Co}\label{co:bdyrevHol}
Let $2<q<\nf$. There exists $\de_0=\de_0(n,m,\La,q)>0$ such that the following holds whenever $\de<\de_0$.
Let $0<R<R_0/8$, and suppose $w\in W^{1,2}(\Om_{2R}(x_0),\R^m)$ solves $L_Aw=0$ with $w=0$ on $\p\Om\cap B_{2R}(x_0)$.
Then there exist positive constants $C_6=C_6(n,m,\La,s)$ and $C=C(n,m,\La)$ such that
\begin{equation}\label{homestLq}
 \bigg(\fint_{\Om_R(x_0)}|Dw|^q\bigg)^{1/q}
 \le C\exp(C_6\mko q^{s-1}\eta^s)
 \bigg(\fint_{\Om_{2R}(x_0)}|Dw|^2\bigg)^{1/2}.
\end{equation}
\end{Co}

\begin{proof}
Set $p_0=2n/3$, $G=|DA|^{p_0}\mkt\mathbf1_{\Om_{2R}(x_0)}$, $p_1\coloneq p_0+1$, and $\alpha_0\coloneq p_0/p_1$.
We denote $\wti w\coloneq w\circ\Phi$.
For $x\in\Om_R(x_0)$, we write $\Psi(x)=(z'(x),z_n(x))$, and set $z_x\coloneq (z'(x),0)$.
By~\eqref{propbiLip} and~\eqref{incimage}, we have $z_x\in\Ga_{R_0/4}$ and
\begin{align}\label{lipgracon}
\begin{dcases}
B_{R/4}^+(z_x)\subset\Psi(\Om_{2R}(x_0)),\\[0.2ex]
    \frac 34 \dist (x,\p\Om)\le z_n(x)\le \frac 54\dist (x,\p\Om)\le \frac 54 R.
\end{dcases}
\end{align}
Since $B_{z_n(x)/4}(x)\subset \Om$, applying the interior estimate~\eqref{Nguyenptbd} and letting $r\to 0$ yields that for almost every $x\in \Om_R(x_0)$,
\begin{equation}\label{1stptbdwbry}
    |Dw(x)|\le C\exp\bigg(C\int_0^{z_n(x)/4} \vt_A(x,t)\,\frac{dt}{t}\bigg)
 \bigg(\fint_{B_{z_n(x)/4}(x)}|Dw|^2\bigg)^{1/2},
\end{equation}
where $C=C(n,m,\La)$.
We now assume $z_n(x)<R/32$. Since $|\Psi(x)-z_x|=z_n(x)$, by~\eqref{propbiLip} and~\eqref{incimage}, we have $\Psi(B_{z_n(x)/4}(x))\subset B_{2z_n(x)}^+(z_x)$. Applying Lemma~\ref{lm:Ngubdy} together with a change of variables gives
\begin{align}\label{ballcompbryest}
\begin{aligned}
\fint_{B_{z_n(x)/4}(x)}|Dw|^2&\le C \fint_{B_{2z_n(x)}^+(z_x)}|D\wti w|^2\\
     &\le C\bigg(\frac R{z_n(x)}\bigg)^{C\de}
 \exp\bigg(C\int_{2z_n(x)}^{R/4}\vt_A(\Phi(z_x),t)\,\frac{dt}{t}\bigg)
 \fint_{B_{R/4}^+(z_x)}|D\wti w|^2.
 \end{aligned}
\end{align}
Moreover, for $2z_n(x)\le t\le R/4$, we have $|x-\Phi(z_x)|\le 4|\Psi(x)-z_x|/3=4z_n(x)/3<t$, and hence $B_{t}(\Phi(z_x) )\subset B_{2t}(x) $. It follows that
\begin{align}\label{vtcompar}
    \vt_A(\Phi(z_x),t)\le C\mko  \vt_A(x,2t),\qquad \forall \, 2z_n(x)\le t\le R/4,
\end{align}
where $C=C(n)$. Since $B_{R/4}^+(z_x)\subset\Psi(\Om_{2R}(x_0))$, combining~\eqref{1stptbdwbry}--\eqref{vtcompar} with a change of variables implies 
\begin{align}\label{finptestDw}
\begin{aligned}
     |Dw(x)|&\le C \bigg(\frac R{z_n(x)}\bigg)^{\wti C\de}\exp\bigg(C\int_0^{R/2} \vt_A(x,t)\,\frac{dt}{t}\bigg)
 \bigg(\fint_{\Om_{2R}(x_0)}|Dw|^2\bigg)^{1/2}\\
 &\le C \bigg(\frac R{z_n(x)}\bigg)^{\wti C\de}\exp\Big(C\mathbf W_{\al_0,p_1}G(x)\Big)\bigg(\fint_{\Om_{2R}(x_0)}|Dw|^2\bigg)^{1/2},
 \end{aligned}
\end{align}
where the constants $C$ and $\wti C$ depend only on $n,m,\La$.
If $z_n(x)\ge R/32$, then by~\eqref{1stptbdwbry}, the estimate~\eqref{finptestDw} remains valid. 

Finally, for $2<q<\nf$, we choose $\de_0=\de_0(n,m,\La,q)$ such that $4\mko\wti C q\mko\de <1$. Then by~\eqref{propbiLip} and~\eqref{lipgracon}, we have 
\[
    \fint_{\Om_R(x_0)}\bigg(\frac R{z_n(x)}\bigg)^{2\mko\wti C q\mko \de}\le C,
\]
where $C=C(n)$. By Cauchy--Schwarz inequality and Lemma~\ref{lm:Wolff-exp}, we obtain
\begin{align*}
&\bigg(\fint_{\Om_R(x_0)}|Dw|^q\bigg)^{1/q}\\
 &\le C  \Bigg(\fint_{\Om_R(x_0)}\bigg(\frac R{z_n(x)}\bigg)^{2\mko\wti C q\mko \de}\Bigg)^{\frac 1{2q}}\bigg(\fint_{B_R(x_0)} \exp \Big(2\mko C\mko q\mko W_{\alpha_0,p_1}G\mko\Big)\bigg)^{\frac 1{2q}}\bigg(\fint_{\Om_{2R}(x_0)}|Dw|^2\bigg)^{1/2}\\
 &\le C\exp(C_6\mko q^{s-1}\eta^s)
 \bigg(\fint_{\Om_{2R}(x_0)}|Dw|^2\bigg)^{1/2},
\end{align*}
where $C=C(n,m,\La)$ and $C_6=C_6(n,m,\La,s)$. This completes the proof.
\end{proof}
\subsection{Proofs of the main results}
\
\vskip5pt
Following the real-variable argument developed in~\cite{CaffaPeral98,Wang03,Shen05,Shen07}, we are now ready to prove our main theorem.
\begin{proof}[Proof of Theorem~\ref{th:divest}]
For $p=2$, the claim follows from~\eqref{eq:W-1pestLA},~\eqref{eq:intcond}, and the Lax--Milgram theorem. 
Let $p\ge 4$, and $f\in L^p(\Om,\R^{m\times n})$.
We denote by $T_A$ the operator on $L^2(\Om,\R^{m\times n})$ that maps $ f$ to $Du$, where $u\in W^{1,2}_0(\Om,\R^m)$ solves~\eqref{sysinth}.  
Then we have
\begin{equation}\label{L2bdT}
    \|T_Af\|_{L^2(\Om)}\le \La\|f\|_{L^2(\Om)}.
\end{equation}
Let $q\in (p,\nf)$. 
We assume $0<\de<\de_0$, where $\de_0=\de_0(n,m,\La,q)$ is defined as in Corollary~\ref{co:bdyrevHol}.
Applying Corollaries~\ref{co:intrevHol} and~\ref{co:bdyrevHol} together with a covering argument, we obtain that, for any $x\in\ov\Om$, $0<r<R_0/32$, and $f$ supported outside $\Om_{2r}(x)$,
\begin{equation}\label{globrevhol}
    \bigg(r^{-n}\int_{\Om_r(x)}|T_A f|^q\bigg)^{1/q}
    \le N\bigg(r^{-n}\int_{\Om_{2r}(x)}|T_Af|^2\bigg)^{1/2},
\end{equation}
where $N= C\exp(C_6\mko q^{s-1}\eta^s)>1$, with positive constants $C_6=C_6(n,m,\La,s)$ and $C=C(n,m,\La)$.
We now follow the proof of~\cite[Thm.~3.3]{Shen05} and keep track of the constants in the estimate. 

We choose a cube $Q_0$ containing $\Om$, and extend $T_Af$ and $f$ by zero outside $\Om$. 
Denote by $2Q_0$ the cube with the same center as $Q_0$ but with twice the side length.
The localized Hardy--Littlewood maximal operator $\mca M_{2Q_0}$ is defined by
\[
    \mca M_{2Q_0}\mko g(x)
    =\sup_{x\in Q\subset 2Q_0}
    \fint_Q |g(y)|\,dy,\qquad \forall\, g\in L^1(2Q_0),\;\forall\,x\in 2Q_0,
\]
where $Q$ ranges over all cubes with sides parallel to the coordinate axes, containing $x$, and contained in $2Q_0$.
For $\la>0$, we set
\[
    E(\la)\coloneq \big\{x\in Q_0:\mca M_{2Q_0}(|T_Af|^2)(x)>\la\big\}.
\] 
 Let $ \tau= (2\mko \vae^{2/p})^{-1}>5^n$, where $\vae\in(0,1)$ is a sufficiently small constant to be chosen later. 
The weak type $(1,1)$ estimate for $\mca M_{2Q_0}$, combined with~\eqref{L2bdT}, implies that
\begin{equation}\label{weak11est}
    |E(\tau\la)|\le \frac{\wti C\La^2}{\tau\la}\int_\Om|f|^2,
    \qquad\forall\,\la>0,
\end{equation}
where $\wti C=\wti C(n)$. Set
\[
    \la_1\coloneq \frac{\wti C\La^2C_{\Om}}{\vae \mko \tau\mko |Q_0|}\int_\Om|f|^2,
\]
where $C_\Om>1$ is a constant depending on $\Om$ and is to be determined.
Let $\la\ge\la_1$. Then by~\eqref{weak11est}, we have  $|E(\tau\la)|\le\vae\mko|Q_0|/C_\Om<\vae\mko |Q_0|$.
Applying the Calder\'on--Zygmund decomposition to $E(\tau\la)$ in $Q_0$, we obtain disjoint open dyadic cubes $\{Q_k\}$ such that for each $k$, it holds that
\begin{equation}\label{Qkcond}
   Q_k\subset Q_0,\qquad  |Q_k\cap E(\tau\la)|>\vae|Q_k|,
    \qquad |\wih Q_k\cap E(\tau\la)|\le\vae|\wih Q_k|,
\end{equation}
where $\wih Q_k$ denotes the dyadic parent of $Q_k$. In particular, this implies 
\[
    |Q_k|<\frac{|E(\tau\la)|}{\vae}\le \frac{|Q_0|}{ C_\Om}.
\]
Thus, we may choose $C_\Om$ large enough such that 
\begin{equation}\label{Qkbd}
    4\diam(Q_k)< \min \Big\{\frac{R_0}{32},\frac{\diam(Q_0)}{5\sqrt n}\Big\}.
\end{equation}
Let $\ga>0$ be a sufficiently small constant to be chosen later. 
We claim that if $\wih Q_k$ contains a point where
$\mca M(|f|^2)\le\ga\la$, then
$\wih Q_k\subset E(\la)$.
Suppose otherwise, then there exist
$x_1,x_2\in\wih Q_k$ such that
\[
    \mca M_{2Q_0}(|T_Af|^2)(x_1)\le\la,
    \qquad
    \mca M(|f|^2)(x_2)\le\ga\la,
\]
where $\mca M$ is the Hardy--Littlewood maximal operator defined by cubes.
By~\eqref{Qkcond}, we can pick a point $y\in Q_k\cap E(\tau\la)$.
By the definition of $E(\tau\la)$, there exists
a cube $Q\subset 2Q_0$ containing $y$ such that
\begin{equation}\label{Tf2lowbd}
    \fint_Q |T_Af|^2>\tau\la.
\end{equation}
If $|Q|\ge |Q_k|$, then since $x_1,y\in \wih Q_k$,
we may enlarge $Q$ to a cube $\wti Q$ satisfying $Q\cup \wih Q_k\subset \wti Q\subset 2Q_0$ and $|\wti Q|\le 3^n|Q|$. 
This contradicts the assumption
$x_1\notin E(\la)$, since $\tau>3^n$.
Hence, we have $|Q|<|Q_k|$. 
Since $T_Af$ vanishes outside $\Om$,
by~\eqref{Tf2lowbd}, we obtain
\begin{equation}\label{QkOmclose}
   \dist (Q_k,\Om)
    \le \diam(Q)\le \diam (Q_k).
\end{equation}
 Set $r_k=4\diam (Q_k)$. By~\eqref{QkOmclose}, we can choose $y_k\in \ov \Om$ such that $Q_k\subset B_{r_k/2}(y_k)$. 
 Then by~\eqref{Qkbd}, we have $\wih Q_k\subset B_{r_k}(y_k)\subset B_{2r_k}(y_k)\subset 2Q_0$.
Set
\[
    f_1=f\mko \mathbf 1_{\Om_{2r_k}(y_k)},
    \qquad f_2=f-f_1.
\]
By the choice of $x_1,x_2$, we have
\begin{equation}\label{Tff1est}
    \int_{\Om_{2r_k}(y_k)}|T_Af|^2
    \le C\la|Q_k|,
    \qquad
    \int_{\Om}|f_1|^2
    \le C\ga\mko\la|Q_k|,
\end{equation}
where $C=C(n)$.
The estimate~\eqref{L2bdT} then gives
\begin{equation}\label{Tf1L2bd}
    \int_\Om|T_Af_1|^2
    \le C\La^2\ga\mko\la|Q_k|.
\end{equation}
Hence by~\eqref{Tff1est}, we obtain
\[
    \int_{\Om_{2r_k}(y_k)}|T_Af_2|^2
    \le C(1+\La^2\ga)\la|Q_k|.
\]
Since $f_2$ vanishes on $\Om_{2r_k}(y_k)$,
the reverse H\"older inequality~\eqref{globrevhol} together with~\eqref{Qkbd} yields
\begin{equation}\label{Tf2Lqbd}
    \int_{\Om_{r_k}(y_k)}|T_Af_2|^q
    \le (CN)^q(1+\La^2\ga)^{q/2}
    \la^{q/2}|Q_k|.
\end{equation}
Moreover, the discussion above implies that for any cube $Q\subset 2Q_0$ intersecting $Q_k\cap E(\tau\la)$ with~\eqref{Tf2lowbd}, it holds that $|Q|\le |Q_k|$ and hence $Q\subset B_{r_k}(y_k)$.
It follows that 
\begin{equation}\label{Etaulaloclar}
    Q_k\cap E(\tau\la)\subset \big\{x\in Q_k:\mca M\big(|T_Af|^2\mathbf 1_{B_{r_k}(y_k)}\big)(x)>\tau\la\big\}.
\end{equation}
Combining~\eqref{Tf1L2bd}--\eqref{Etaulaloclar} with the weak type $(1,1)$ and $(q/2,q/2)$ estimates for the maximal operator, we obtain
\begin{equation}\label{eq:Shenlocal}
    |Q_k\cap E(\tau\la)|
    \le \wti C|Q_k|\bigg(
    \frac{\La^2\ga}{\tau}
    +\frac{(\wti CN)^q(1+\La^2\ga)^{q/2}}
    {\tau^{q/2}}\bigg),
\end{equation}
where $\wti C=\wti C(n)>0$.
We choose
\begin{equation}\label{defvaega}
    \vae=\Big(\frac {c_0}{N}\Big)^{\frac{pq}{q-p}},\qquad
    \ga={c_0}\La^{-2}\tau\mko\vae,
\end{equation}
where $c_0=c_0(n)\in (0,1)$ is sufficiently small. In particular, we have 
\[
    \tau=(2\mko \vae^{2/p})^{-1}=2^{-1} \Big(\frac {c_0}{N}\Big)^{-\frac {2q}{q-p}}>(2c_0^2)^{-1}.
\]
Thus, the condition $\tau>(10\sqrt n)^n$ is satisfied for sufficiently small $c_0$. If $c_0<(2\wti C)^{-1}$, then we have
\begin{equation}\label{CN12gatau}
    \wti C \frac{\La^2\ga}{\tau}<\frac 12\mko \vae.
\end{equation}
In addition, since $p\ge 4$, it holds that \begin{equation}\label{N12ga<1}
    \La^2\ga=c_0\tau\mko\vae<c_0\mko \vae^{1/2}<1.
\end{equation}
The definition~\eqref{defvaega} of $\vae$ gives $N=c_0 \mko\vae ^{1/q-1/p}$. It follows that
\[
    \wti C N \Big(\frac{1+\La^2 \ga}{\tau}\Big)^{1/2}\le  2\mko\wti C N\mko \vae^{1/p}=2\mko\wti C \mko c_0\mko \vae^{1/q}.
\]
Since $q>p\ge 4$, for $c_0=c_0(n)$ small enough, we then obtain
\begin{equation}\label{CN2N1taubd}
   \wti C \frac{(\wti CN)^q(1+\La^2\ga)^{q/2}}
    {\tau^{q/2}}\le \wti C (2\wti C c_0)^q \mko \vae\le \wti C (2\wti C c_0)^4 \mko \vae<\frac \vae 2.
\end{equation}
Combining~\eqref{CN12gatau} and~\eqref{CN2N1taubd} with~\eqref{eq:Shenlocal} yields $|Q_k\cap E(\tau\la)|\le\vae|Q_k|$, which contradicts~\eqref{Qkcond}.
Thus, any $\wih Q_k$ intersecting the set $\{\mca M(|f|^2)\le\ga\la\}$ is contained in $E(\la)$. 
Since any two dyadic cubes are either disjoint or one is contained in the other, by taking a disjoint subcover $\{\wih Q_{k'}\}$ of $E(\tau\la)\cap \{\mca M(|f|^2)\le\ga\la\}$, we obtain that
\begin{equation}\label{EtaulaMf2com}
    |E(\tau\la)|\le\vae|E(\la)|
    +\big|\{\mca M(|f|^2)>\ga\la\}\big|,
    \qquad\forall\,\la\ge\la_1.
\end{equation}
Since $p\ge 4$, we have $\|\mca M\|_{L^{p/2}(\R^n)\to L^{p/2}(\R^n)}\le C$, where $C$ is independent of $p$. It follows that for any $\mu  \ge \tau\la_1$,
\begin{align*}
   & \frac p2(1-\vae \mko\tau^{p/2})\int_0^{\mu}\la^{p/2-1}|E( \la)|\,d\la\\
    &\le \frac p2\mko|Q_0|\int_0^{\tau\la_1}\la^{p/2-1}\,d\la+\frac p2 \int_0^{\nf}\la^{p/2-1}\Big|\Big\{\mca M(|f|^2)>\frac{\ga\la}\tau\Big\}\Big|\,d\la\\
    &=(\tau\la_1)^{p/2}|Q_0|+\Big(\frac\tau \ga\Big)^{p/2}\int_{\R^n} \big|\mca M(|f|^2)\big|^{p/2}\\
    &\le (\tau\la_1)^{p/2}|Q_0|+\Big(\frac{C\tau} \ga\Big)^{p/2}\int_{\Om} |f|^p,
\end{align*}
where $C=C(n)$. 
Since $\vae \mko \tau^{p/2}=2^{-p/2}\le 1/4$, letting $\mu\to \nf$ yields $T_Af\in L^p(\Om)$ with 
\begin{align*}
    \int_{\Om} |T_Af|^p&\le \int_{Q_0}  |\mca M_{2Q_0}(|T_Af|^2)|^{p/2}=\frac p2\int_0^{\nf}\la^{p/2-1}|E( \la)|\,d\la\\
    &\le \frac 43 (\tau\la_1)^{p/2}|Q_0|+\frac 43\Big(\frac{C\tau} \ga\Big)^{p/2}\int_{\Om} |f|^p.
\end{align*}
By the definitions of $\la_1,\tau,\ga$, and H\"older's inequality, we then obtain
\begin{align}\label{TfLpbd}
\begin{aligned}
    \|T_Af\|_{L^p(\Om)}
    &\le C_{\Om}'\mkt \La\mko \vae^{-1/2}\|f\|_{L^p(\Om)}\\
    &\le C_\Om'\mkt \La
    \Big(\frac N{c_0}\Big)^{\frac{pq}{2(q-p)}}\|f\|_{L^p(\Om)},
    \end{aligned}
\end{align}
where $C_\Om'$ depends only on $\Om$.
Finally, we choose $q=2p$. Combining~\eqref{TfLpbd} with the expression $N= C\exp(C_6\mko q^{s-1}\eta^s)$ in~\eqref{globrevhol} yields the required estimate for $p\ge 4$.
The cases $2<p<4$ and $1<p<2$ follow by interpolation and duality, respectively.
\end{proof}
\begin{proof}[Proof of Theorem~\ref{Th:divHajespa}]
We adapt the proof of Theorem~\ref{th:divest}. Set
\[
    \wti p_0=\frac{2n}{1+2\theta},\qquad
    \wti p_1=\wti p_0+1,\qquad
    \wti\al=\frac{\theta\wti p_0}{\wti p_1}.
\]
For $x\in\ov\Om$ and $r>0$, we define
\[
    \wti\vt_A(x,r)\coloneq r^\theta
    \bigg(\fint_{\Om_r(x)}\wti g^{\wti p_0}\bigg)^{1/\wti p_0}.
\]
Since $1/\wti p_0-\theta/n=1/(2n)$, the fractional Hajłasz--Sobolev--Poincaré inequality (see~\cite{Zhou11,Bojarski87,Hajłasz03}) implies that for all $x\in \ov \Om$ and $0<r<R_0/8$,
\begin{align*}
    \bigg(\fint_{\Om_r(x)}
    |A-(A)_{\Om_r(x)}|^{2n}\bigg)^{1/(2n)}
    \le C\mko\wti\vt_A(x,r),
\end{align*}
where $C=C(n,\theta)$. 
Thus, after subtracting a constant spatially antisymmetric matrix, the coefficient estimates~\eqref{L2nnormA},~\eqref{Poinoscder},~\eqref{eq:bdycoeffL2n}, and~\eqref{eq:comcoeffbd} remain valid with $\vt_A$ replaced by $\wti\vt_A$, where the constants depend additionally on $\theta$. 
Consequently, the proofs of Lemmas~\ref{lm:L2grad},~\ref{lm:intcom}, and~\ref{lm:Nguint}, together with their boundary analogues Lemmas~\ref{L2firgrabdy},~\ref{lm:bdycom}, and~\ref{lm:Ngubdy}, apply with $\wti\vt_A$ in place of $\vt_A$.

Now let $x_0\in \ov \Om$, $0<R<R_0/8$, and set $G=\wti g^{\wti p_0}\mathbf1_{\Om_{2R}(x_0)}$. 
For any $x\in \Om_R(x_0)$, we have 
\[
    \int_0^{R}\wti\vt_A(x,t)\,\frac{dt}{t}
    \le C\mko\mathbf W_{\wti\al,\wti p_1}G(x),
\]
where $C=C(n,\theta)$. 
We now apply Lemma~\ref{lm:Wolff-exp} with $\al=\wti\al$ and $p=\wti p_1$. Since
\[
    \wti\al\wti p_1=\theta\wti p_0<n,
    \qquad
    \frac{n(\wti p_1-1)}{\wti\al\wti p_1}=\frac n\theta,
\]
we obtain that for any $\ga>0$,
\begin{align*}
    \fint_B\exp\big(\ga\mathbf W_{\wti\al,\wti p_1}
   G\big)
    \le4\exp\big(C(\ga\eta)^s\big),
\end{align*}
where $C=C(n,\theta,s)$. 
The proofs of Corollaries~\ref{co:intrevHol} and~\ref{co:bdyrevHol} then yield the same reverse Hölder estimates as in~\eqref{globrevhol} with
\[
    N=C\exp(C_7\mko q^{s-1}\eta^s),
\]
where $C=C(n,m,\La)$ and $C_7=C_7(n,m,\La,\theta,s)$. 
The rest of the proof follows from the same real-variable argument as in~\cite{Shen05} and the proof of Theorem~\ref{th:divest}.
\end{proof}
\section{Sharpness of the estimate~\eqref{W1pestW1n}}\label{sec:exa}
In this section, we show that the exponential dependence on $\eta^s$ in~\eqref{W1pestW1n} is sharp. 
For scalar equations, the construction is based on the radial examples of Meyers~\cite{Meyers63} and Serrin~\cite{Serrin64}, and the computation in~\cite{JinMazyaVanSch09}.
For systems, we use a parameterized version of De Giorgi's example; see also~\cite[Ch.~II, Ex.~3.1]{Giaquinta83} and~\cite[Eq.~(2.16)]{Soucek84}.
We truncate these $W^{1,(n,\nf)}$ coefficients near the origin to obtain sequences of smooth coefficients satisfying the assumptions of Theorem~\ref{th:divest}.

For a coefficient tensor $A$ on $B_1$, we denote by $T_A$ the operator mapping $f$ to $Du$, where $u\in W^{1,2}_0(B_1,\R^m)$ solves
\[
    \p_\al\big(A^{\al\beta}_{ij}\p_\beta u^j\big)
    =\p_\al f_i^\al\qquad\text{in }B_1.
\]
For scalar equations, we prove that the quantitative dependence on $\eta$ is sharp for $p>n$.
\begin{Prop}\label{prop:shasca}
Let $2\le n<p<\nf$ and $1\le s<\nf$. 
Then there exist constants $\La=\La(n,p)\ge1$, $\vae_0=\vae_0(n,p)>0$, and symmetric coefficient matrices $A_\vae\in C^\nf(\ov{B_1},\R^{n\times n})$ for all $0<\vae<\vae_0$, such that
\begin{equation}\label{Aepell}
    \La^{-1}|\xi|^2\le
    A_\vae^{\al\beta}(x)\xi_\beta\xi_\al
    \le\La|\xi|^2,
    \qquad \forall\,x\in B_1,\quad \forall\,\xi\in\R^n,
\end{equation}
and, writing $\eta_\vae=    \|DA_\vae\|_{L^{n,s}(B_1)}$, we have
\begin{gather}
\begin{gathered}\label{scaLpgrow}
    c\log\frac1\vae\le\eta_\vae^s\le C\log\frac1\vae,\\
    \|T_{A_\vae}\|_{L^p(B_1)\to L^p(B_1)}
    \ge c\mko \exp(c\mkt\eta_\vae^s),
\end{gathered}
\end{gather}
where $c,C>0$ depend only on $n,p,s$.
\end{Prop}

\begin{proof}
Fix $\si\in (n/p,1)$. For $x\ne0$, we set
\begin{equation}\label{scaconA}
    A^{\al\beta}(x)=\de_{\al\beta}-\frac{\si(n-\si)}{n-1}
    \Big(\de_{\al\beta}-\frac{x_{\al} x_\beta}{|x|^2}\Big),
    \qquad w(x)=x_1|x|^{-\si}.
\end{equation}
When $n=2$, this construction coincides with that in~\cite[Eqs.~(65)--(66)]{Meyers63}.
The eigenvalues of $A$ are $1$ with multiplicity $1$, and
\[
    1-\frac{\si(n-\si)}{n-1}
    =\frac{(1-\si)(n-1-\si)}{n-1}>0
\]
with multiplicity $n-1$.
Thus $A$ is uniformly elliptic. 
Applying~\cite[Lem.~2.1]{JinMazyaVanSch09} with $v(r)=r^{-\si}$ and $\al(r)=-\si(n-\si)/(n-1)$ gives
\[
    \p_\al(A^{\al\beta}\p_\beta w)=0
    \qquad\text{in }\;B_1\smi\{0\}.
\]
Since $\si< 1\le n/2$, we have $w\in W^{1,2}(B_1)$. 
By~\cite[Ch.~II, Lem.~3.1]{Giaquinta83}, it follows that
\begin{equation}\label{wsoldis}
     \p_\al(A^{\al\beta}\p_\beta w)=0
    \qquad\text{in }\;\mca D'(B_1).
\end{equation}
Let $\zeta\in C_c^\nf(B_{3/4})$ with $\zeta=1$ on $B_{1/2}$, and set $u=\zeta w$. Then $L_Au\in C_c^\nf(B_1)$ is supported in $B_1\smi B_{1/2}$.
We choose $f=\{f_\al\}\in C^\nf(\ov{B_1},\R^n)$ such that
$L_Au=\p_\al f_\al$.
Let $\chi\in C^\nf([0,\nf))$ satisfy $0\le\chi\le1$, $\chi=0$ on $[0,1]$, and $\chi=1$ on $[2,\nf)$. For $0<\vae<1/8$, define
\begin{equation}\label{Aepsdef}
    A_\vae(x)=I_n+\chi(|x|/\vae)\big(A(x)-I_n\big).
\end{equation}
Then each $A_\vae$ satisfies~\eqref{Aepell} with a fixed $\La=\La(n,p)\ge 1$. Let $u_\vae\in W^{1,2}_0(B_1)$ solve $L_{A_\vae}u_\vae=\p_\al f_\al$. 
Since $A_\vae$ and $f_\al$ are smooth, we have $u_\vae\in C^\nf(\ov{B_1})$. Subtracting the equations gives
\[
    L_{A_\vae}(u_\vae-u)
    =\p_\al\big(\mko(A^{\al\beta}-A_\vae^{\al\beta})\p_\beta u\big).
\]
Since $A_\vae=A$ outside $B_{2\vae}$, the $L^2$ estimate implies
\begin{equation}\label{uvaeudifest}
    \|D(u_\vae-u)\|_{L^2(B_1)}
    \le C\|Du\|_{L^2(B_{2\vae})}
    \le \wti C_0\mko\vae^{n/2-\si},
\end{equation}
where $\wti C_0=\wti C_0(n,p)>0$.
On the other hand, homogeneity gives
\[
    \|Du\|_{L^2(B_\rho\smi B_{\rho/2})}
    =c_0\mko\rho^{n/2-\si},\qquad \forall\, 0<\rho<1/2,
\]
where $c_0=c_0(n,p)>0$. 
We fix $R=R(n,p)>4$ sufficiently large such that $c_0R^{n/2-\si}\ge2\wti C_0$. 
For $\vae<(2R)^{-1}$, we obtain
\[
    \|Du_\vae\|_{L^2(B_{R\vae}\smi B_{R\vae/2})}
    \ge  \|Du\|_{L^2(B_{R\vae}\smi B_{R\vae/2})}-   \|D(u_\vae-u)\|_{L^2(B_1)}\ge \frac{c_0}{2}(R\vae)^{n/2-\si}.
\]
H\"older's inequality then implies
\begin{align*}
    \|Du_\vae\|_{L^p(B_1)}
    \ge \|Du_\vae\|_{L^p(B_{R\vae}\smi B_{R\vae/2})}  \ge c_1\mko \vae^{n/p-\si},
\end{align*}
where $c_1=c_1(n,p)>0$.
Since $f$ is fixed, this proves
\begin{equation}\label{Aepnormest}
    \|T_{A_\vae}\|_{L^p(B_1)\to L^p(B_1)}
    \ge c_2\mkt \vae^{-(\si-n/p)},
\end{equation}
where $c_2=c_2(n,p)>0$.

It remains to estimate $\eta_\vae$. By~\eqref{Aepsdef}, we have $|DA_\vae(x)|\le C/\max(\vae,|x|)$, where $C=C(n,p)$.
The definition~\eqref{eq:defLor} of Lorentz spaces then gives
\begin{equation}\label{etaepupbd}
    \eta_\vae^s
    \le C\vae^{-s}\int_0^{|B_1|\vae^n}t^{s/n}\,\frac{dt}{t}
       +C\int_{|B_1|\vae^n}^{|B_1|}\frac{dt}{t}
    \le C\Big(1+\log\frac1\vae\Big),
\end{equation}
where $C=C(n,p,s)$.
Since $A_\vae=A$ on $B_1\smi B_{2\vae}$, by homogeneity, we obtain positive constants $\wti C_1,c_3,c_4$ depending only on $n,p$ such that  
\begin{equation}\label{lowbdAep}
     d_{|DA_\vae|}(\la)\ge c_4  \la^{-n},\qquad \forall\, \wti C_1\le \la \le c_3\mkt \vae^{-1}.
\end{equation}
It follows that $\eta_\vae^s\ge c_5\log(1/\vae)$ for small enough $\vae$, where $c_5=c_5(n,p,s)>0$. 
Combined with~\eqref{Aepnormest}, this proves~\eqref{scaLpgrow}.
\end{proof}
For elliptic systems, the quantitative dependence on $\eta$ is sharp for any $p\neq 2$.
By duality, it suffices to consider $p>2$.
\begin{Prop}\label{prop:sharpsys}
Let $n\ge2$, $2<p<\nf$, and $1\le s<\nf$.
Then there exist constants $\La=\La(n,p)\ge1$, $\vae_0=\vae_0(n,p)>0$, and coefficient tensors $A_\vae\in C^\nf(\ov{B_1},\R^{n\times n\times n\times n})$ for all $0<\vae<\vae_0$, satisfying $(A_\vae)^{\al\beta}_{ij}=(A_\vae)^{\beta\al}_{ji}$ and
\begin{equation}\label{LegAep}
    \La^{-1}|\xi|^2\le
    (A_\vae)^{\al\beta}_{ij}(x)\xi_\beta^j\xi_\al^i
    \le\La|\xi|^2,
    \qquad \forall\, x\in B_1,\quad \forall\,\xi\in\R^{n\times n}.
\end{equation}
In addition, writing  $\eta_\vae=    \|DA_\vae\|_{L^{n,s}(B_1)}$, we have
\begin{gather}\begin{gathered}\label{sysLpgrow}
    c\log\frac1\vae\le\eta_\vae^s\le C\log\frac1\vae,\\
    \|T_{A_\vae}\|_{L^p(B_1)\to L^p(B_1)}
    \ge c\mko\exp(c\mko\eta_\vae^s),
    \end{gathered}
\end{gather}
where $c,C>0$ depend only on $n,p,s$.
\end{Prop}

\begin{proof}
Fix $\si\in (n/p,n/2)$.
For $x\ne0$, we set
\begin{equation}\label{exsysA}
\begin{aligned}
    B_i^\al(x)&=(n-2)\delta_{i\al}+n\frac{x_i x_\al}{|x|^2},&\qquad\kappa&=\frac{\si(n-\si)}{(n-1)^2(n-2\si)^2},\\
    A^{\al\beta}_{ij}(x)&=\delta_{ij}\delta_{\al\beta}+\kappa B_i^\al(x)B_j^\beta(x),&\qquad
    w^i(x)&=x_i\mko |x|^{-\si}.
\end{aligned}
\end{equation}
When $\kappa=1$ and $n\ge 3$, this coincides with De Giorgi's original construction~\cite{DeGiorgi68}. 
Since $|B|^2=n^2(n-1)$, we have
\[
    |\xi|^2\le A^{\al\beta}_{ij}\xi_\beta^j\xi_\al^i
    =|\xi|^2+\kappa\big(B_i^\al\xi_\al^i\big)^2
    \le\big(1+\kappa \mko n^2(n-1)\mko\big)|\xi|^2.
\]
Writing $r=|x|$, direct computation yields that on $B_1\smi\{0\}$,
\begin{align*}
    B_j^\beta\p_\beta w^j
    &=(n-1)(n-2\si)r^{-\si},\\
    \p_\al(r^{-\si}B_i^\al)
    &=(n-1)(n-2\si)x_i \mko r^{-\si-2},\\
    \lap w^i&=-\si(n-\si)x_i\mko r^{-\si-2}.
\end{align*}
It follows that on $B_1\smi\{0\}$,
\[
    \p_\al\big(A^{\al\beta}_{ij}\p_\beta w^j\big)
    =\big(-\si(n-\si)+\kappa(n-1)^2(n-2\si)^2\big)
      x_i\mko r^{-\si-2}=0.
\]
Since $\si<n/2$, we have $w\in W^{1,2}(B_1,\R^n)$, and hence by~\cite[Ch.~II, Lem.~3.1]{Giaquinta83}, we obtain
\[
      \p_\al\big(A^{\al\beta}_{ij}\p_\beta w^j\big)=0
    \qquad\text{in }\;\mca D'(B_1,\R^n).
\]
We use the cut-off functions $\zeta$ and $\chi$ from the proof of Proposition~\ref{prop:shasca}, set $u=\zeta w$, and define
\[
    (A_\vae)^{\al\beta}_{ij}(x)
    =\delta_{ij}\delta_{\al\beta}
    +\chi(|x|/\vae)
      \big(A^{\al\beta}_{ij}(x)-\delta_{ij}\delta_{\al\beta}\big).
\]
Each $A_\vae$ satisfies~\eqref{LegAep} with a fixed $\La=\La(n,p)\ge 1$.
Fix $f=\{f_i^\al\}\in C^\nf(\ov{B_1},\R^{n\times n})$ satisfying $L_Au=\p_\al f_i^\al$,
and let $u_\vae\in W^{1,2}_0(B_1,\R^n)$ solve $L_{A_\vae}u_\vae=\p_\al f_i^\al$.
The same argument as in~\eqref{uvaeudifest}--\eqref{lowbdAep} then proves the estimates~\eqref{sysLpgrow}.
\end{proof}

\end{document}